\documentclass[11pt]{article}
\usepackage[english]{babel}
\usepackage{fullpage}

\usepackage{amsmath,amssymb,dsfont,amsthm,cases,bm}
\usepackage{upgreek,fourier-orns} 
\usepackage{graphicx,float,color,subfigure,eso-pic}
\definecolor{urlcol}{rgb}{0.3,0.2,0.3}
\definecolor{linkcol}{rgb}{0.2,0.2,0.5}
\definecolor{citecol}{rgb}{0.1 ,0.5 ,0.1}
\usepackage{hyperref}
\hypersetup{
    pdftoolbar=true,        
    pdfmenubar=true,        
    pdffitwindow=true,      
    pdftitle={CV},    
    pdfauthor={Stéphane Gerbi},     
    pdfsubject={CV},   
    pdfcreator={TeXShop},   
    pdfproducer={TeXShop}, 
    pdfkeywords={Resume of Stéphane Gerbi}, 
    colorlinks=true,       
    linkcolor=blue,  
    citecolor=blue,        
    filecolor=magenta,      
    urlcolor=cyan           
}

\usepackage{url}
\newtheorem{thm}{Theorem}[section]
\newtheorem{lem}{Lemme}[section]

\newtheorem{rque}{Remark}[section]

\def\dsp{\displaystyle}

\def\abs#1{\left\vert #1 \right\vert}

\def\dive{{\rm div}}

\newcommand{\fig}{F{\sc{ig}}. }
\newcommand{\R}{\mathbb{R}}

\newcommand{\vecdeux}[2]{\left(\begin{array}{c} #1\\#2 \end{array}\right)}
\newcommand{\vectrois}[3]{\left(\begin{array}{c}#1\\#2\\#3\end{array}\right)}
\newcommand{\PFS}{\textbf{PFS}}
\newcommand{\FS}{\textbf{FS}}
\newcommand{\Pres}{\textbf{P}}

\newcommand{\FF}{\textbf{F}}
\newcommand{\TT}{\textbf{T}}
\newcommand{\NN}{\textbf{N}}
\newcommand{\BB}{\textbf{B}}

\newcommand{\UU}{{\textbf{U}}}

\newcommand{\ii}{{\textbf i}}
\newcommand{\ji}{{\textbf j}}
\newcommand{\ki}{{\textbf k}}
\newcommand{\nn}{{\textbf n}}
\newcommand{\mm}{{\textbf m}}

\newcommand{\VV}{\textbf{V}}
\newcommand{\SE}{\mathcal{\textbf{S}}}

\newenvironment{@abssec}[1]{%
     \if@twocolumn
       \section*{#1}%
     \else
       \vspace{.05in}\footnotesize
       \parindent .2in
         {\upshape\bfseries #1. }\ignorespaces
     \fi}
     {\if@twocolumn\else\par\vspace{.1in}\fi}
\renewenvironment{abstract}{\begin{@abssec}{Abstract}}{\end{@abssec}}
\newenvironment{keywords}{\begin{@abssec}{Key words}}{\end{@abssec}}
\everymath{\dsp}

\date{}
\title{A mathematical model for unsteady mixed flows in closed water pipes}
\author{Christian Bourdarias\footnote{Laboratoire de Math\'ematiques, Universit\'e de Savoie,  73376 Le Bourget du Lac, France, (Christian.Bourdarias@univ-savoie.fr).}
\and
Mehmet Ersoy\thanks{BCAM - Basque Center for Applied Mathematics, Bizkaia Technology Park, Building 500. 48160 Derio Basque Country - Spain, (mersoy@bcamath.org).}
\and
St\'{e}phane Gerbi\thanks{Laboratoire de Math\'ematiques, Universit\'e de Savoie, 73376 Le Bourget du Lac, France, (Stephane.Gerbi@univ-savoie.fr).}
}

\begin{document}

\maketitle

\begin{abstract}
We present the formal derivation of a new unidirectional model for unsteady mixed flows in non uniform closed water pipe.
In the case of free surface incompressible  flows, the \FS-model is formally obtained, using formal asymptotic analysis, which is an extension to more classical shallow water models. In the same way, when the
pipe is full, we propose the \Pres-model, which describes the evolution of
a compressible inviscid flow, close to gas dynamics equations in a nozzle. In order to cope the transition between a free surface state and a pressured (i.e. compressible) state, we propose a mixed model, the \PFS-model, taking into account changes of section and slope variation.
\end{abstract}
\begin{keywords}
Shallow water equations, unsteady mixed flows, free surface flows, pressurized flows,
curvilinear transformation, asymptotic analysis.
\end{keywords}

{\scriptsize
\subsubsection*{Notations concerning geometrical variables}

\begin{itemize}
\item $(0,\ii,\ji,\ki)$: Cartesian reference frame
\item $\omega(x,0,b(x))$: parametrization in the reference frame $(0,\ii,\ji,\ki)$ of the plane curve $\mathcal{C}$ which corresponds to the main flow axis
\item $(\TT,\NN,\BB)$: Serret-Frenet reference frame attached to $\mathcal{C}$ with $\TT$ the tangent, $\NN$ the normal and $\BB$ the bi-normal vector
\item $X,Y,Z$: local variable in the Serret Frenet reference frame with $X$ the curvilinear abscissa, $Y$ the width of pipe, $Z$ the
$\BB$-coordinate of any particle.
\item $\sigma(X,Z)=\beta(X,Z)-\alpha(X,Z)$: width of the pipe at $Z$ with $\beta(X,Z)$ (resp. $\alpha(X,Z)$) is the right (resp. left) boundary point at altitude $Z$
\item $\theta(X)$:  angle $(\ii,\TT)$
\item $S(X)$:  cross-section area
\item $R(X)$: radius of the cross-section $S(X)$
\item $\nn_{\textbf{wb}}$: outward normal vector to the wet part of the pipe
\item $\nn$: outward normal vector at the boundary point $m$ in the  $\Omega$-plane
\end{itemize}
\subsubsection*{Notations concerning the free surface (FS) part}
\begin{itemize}
\item $\Omega(t,X)$:  free surface cross section
\item $H(t,X)$: physical water height
\item $h(t,X)$: $Z$-coordinate of the water level
\item $\nn_{\textbf{fs}}$: outward $\BB$-normal vector to the free surface
\item $A$:  wet area
\item $Q$:  discharge
\item $\rho_0$:  density of the water at atmospheric pressure $p_0$
\end{itemize}
\subsubsection*{Notations concerning the pressurised part}
\begin{itemize}
\item $\Omega(X)$:  pressurised cross section
\item $\rho(t,X)$:  density of the water
\item $\beta$:  water compressibility coefficient
\item $c=\frac{1}{\sqrt{\beta\,\rho_0}}$:  sonic speed
\item $A = \frac{\rho}{\rho_0} S$:  FS equivalent wet area
\item $Q$:  FS equivalent discharge
\end{itemize}
\subsubsection*{Notations concerning the \textbf{PFS} model}
\begin{itemize}
\item $\SE$: the physical wet area: $\SE=A$ if the state is free surface, $S$ otherwise
\item $\mathcal{H}$: the $Z$ coordinate of the water level: $\mathcal{H}=h$ if the state is free surface, $R$ otherwise
\end{itemize}
\subsubsection*{Other notations}
\begin{itemize}
\item Bold characters are used for vectors except for $\SE$
\end{itemize}
}
%
\section{Introduction}
The presented work takes place in a more general framework: the modelling of unsteady mixed flows
in any kind of closed domain taking into account  the cavitation problem and air entrapment. 
We are
interested in flows occurring in closed pipe with non uniform sections, where some parts of the flow
can be free surface (it means that only a part of the pipe is filled) and other parts are
pressurised (it means that the pipe is full). The transition phenomenon between the two types of
flows occurs in many situations such as storm se\-wers, waste or supply pipes in hydroelectric
installations. It can be induced by sudden change in the boundary conditions as failure pumping.
During this process, the pressure can  reach severe values and  cause damages. 
The simulation of such a phenomenon is thus a major challenge and a great amount of works was devoted to it these last
years (see \cite{N90},\cite{R85},\cite{SWB98},\cite{CSZ97} for instance).
Recently Fuamba \cite{M02} proposed a model for the transition from a free surface flow to a
pressurised one in a way very close to ours.

The classical shallow water equations are commonly used to describe  free surface flows in open channels. They are also used in the study of mixed flows  using the Preissman slot artefact (see for example \cite{CSZ97,SWB98}). However, this technique does not take into account  the  depressurisation phenomenon which occurs during a water hammer. On the other hand  the Allievi equations, commonly used to describe pressurised flows, are written in a non-conservative form which is not well adapted to a natural coupling with the shallow water equations.

A model for the unsteady mixed water flows in closed pipes and a finite volume
discretization have been previously studied by two of the authors \cite{BG07} and a kinetic formulation  has
been proposed in \cite{BGG08}. We propose here the \PFS-model which tends to extend naturally the work in \cite{BG07} in the case of a closed pipe with  non uniform section. For the sake of simplicity, we do not deal with the deformation of the domain induced by the change of pressure. We will consider only an infinitely rigid pipe.

The paper is organized as follows.
The first section is devoted to the derivation of  the free surface model from the 3D incompressible Euler equations which are written in a suitable local reference frame in order to take into account the local effects produced by the changes of section and the slope variation. To this end, we present two models derived by two techniques inspired from the works in \cite{BFML07} and \cite{GP01}. The first one consists in taking the mean value in the Euler equations along the normal section to the main axis. The  obtained model provides a description taking in account the geometry of the domain, namely the changes of section and also the inertia strength produced by the slope variation. The second one is a formal asymptotic analysis. In this approach, we seek for  an approximation at the first order   and, by comparison with the previous model, the term related to the inertia strength vanishes since it is a term of second order. We obtain the FS-model. 
In Section~\ref{SectionDerivationEquationsP}, we follow the derivation of the FS-model and we derive the model for pressurised flows, called P-model, from the 3D compressible Euler equations by a formal asymptotic analysis. Writing the source terms into a unified form and  using the same couple of conservative unknowns as in \cite{BG08}, we propose  in Section~\ref{SectionCoupling}, a natural model  for mixed flows, that we call \PFS-model, which ensures the continuity of the unknowns and the source terms.

\section{Formal derivation of the free surface model}\label{SectionFormalDerivationFSModel}
The classical shallow water equations are commonly  used to describe physical situations like rivers, coastal domains, oceans and sedimentation problems.
These equations are obtained from  the  incompressible  Euler system (see e.g. \cite{AL08,LOE96}) or from the incompressible Navier-Stokes system (see for instance \cite{BCNV08,BN07,GP01,M07}) by several techniques
(e.g. by direct integration or asymptotic analysis or as in \cite{DBMS08} and especially as proposed by Bouchut \emph{et al.} \cite{BFML07,BMCPV03} from which the \PFS-model is based).

In order to formally derive a unidirectional  shallow water type equation for free surface flow in closed water pipe with varying slope and section, we consider
 that the length of the pipe is larger than the diameter and we write the incompressible Euler equations in a local  Serret-Frenet frame attached to a given plane curve (generally the main pipe axis, see
Remark \ref{RemarkRestrictionPlaneCurve}).
Then, taking advantage of characteristic scales, we  perform a thin layer asymptotic analysis with respect to some small  parameter $\varepsilon = \frac{H}{L}$ which is also assumed to be proportional to the vertical, $W$ and horizontal, $U$ ratio of the fluid movements, i.e. $\varepsilon = \frac{W}{U}$. This assumption translates the fact that in such domain, the flow follows  a main flow axis.
Finally, the equations are vertically averaged  along orthogonal sections to the given plane curve and  we get the \textbf{F}ree \textbf{S}urface model called \FS-model.

Throughout this section, we only consider pipes with variable circular section. However, this analysis can be easily adapted to any type of closed pipes.

\begin{figure}[H]
 \begin{center}
 \includegraphics[scale = 0.6]{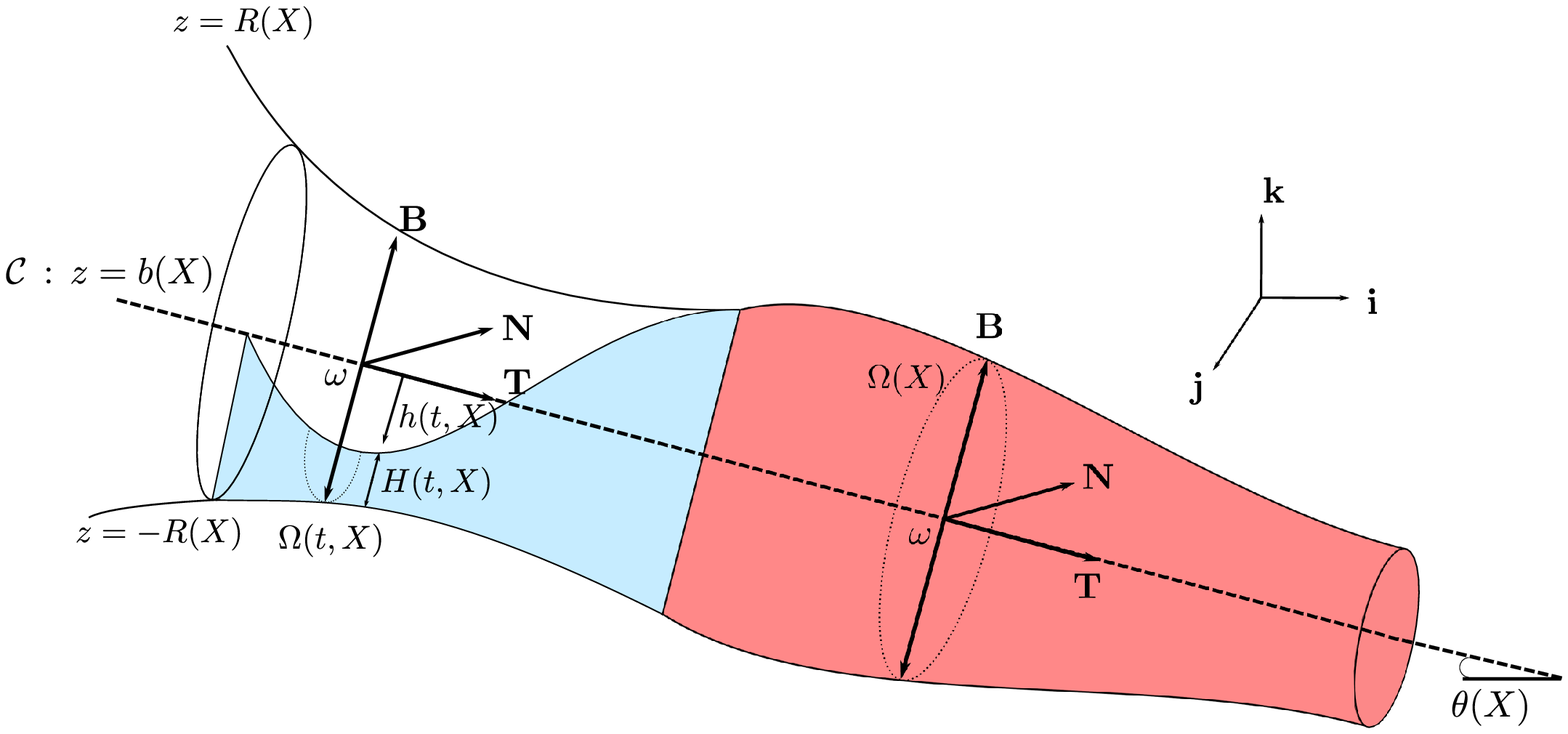}
  \caption{Geometrics characteristics of the pipe}
  \small{Mixed flows: free surface and pressurized}
  \label{FigureOxOz}
 \end{center}
\end{figure}

\begin{figure}[H]
 \begin{center}
 \includegraphics[scale = 0.45]{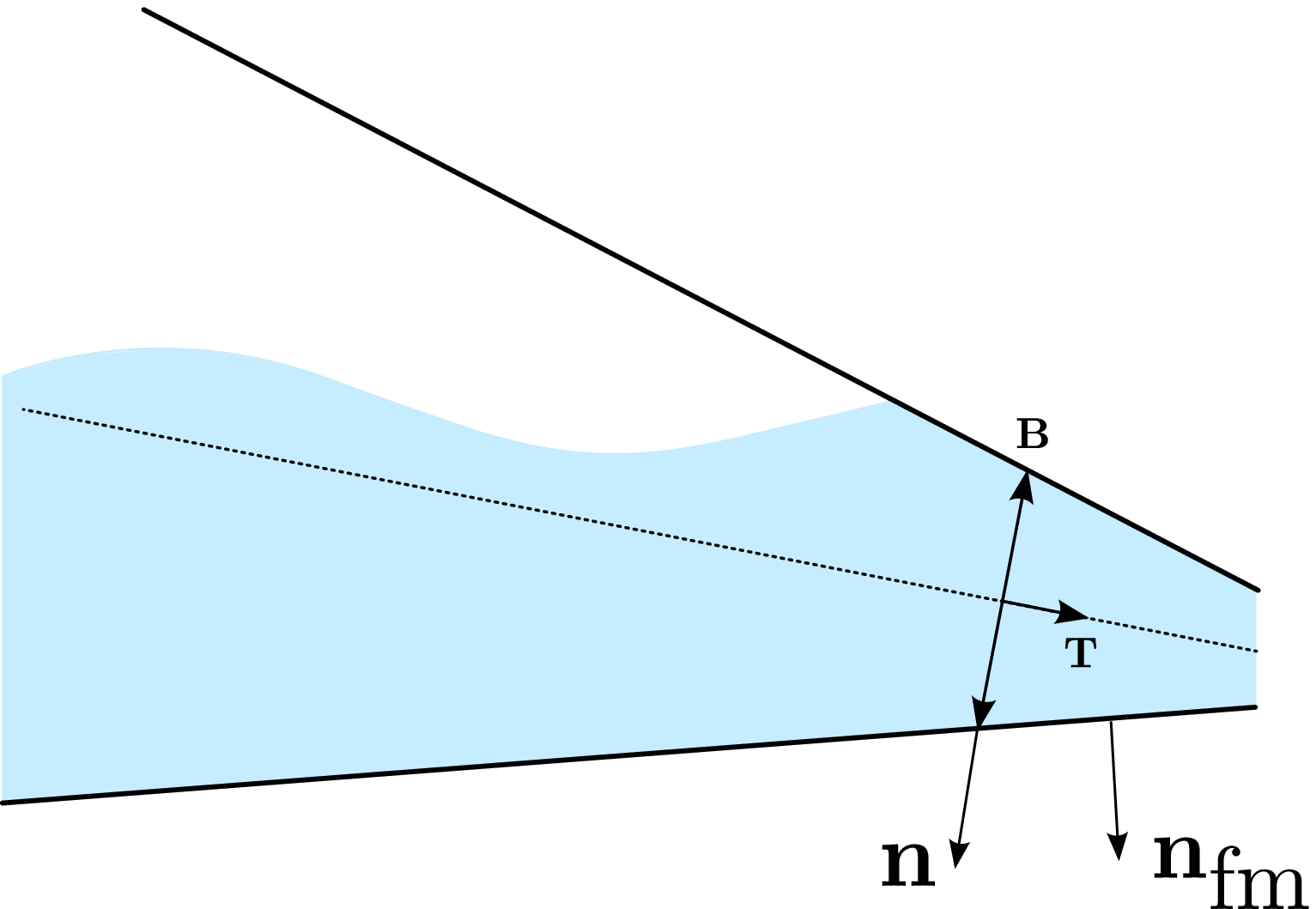}
  \caption{outward unit normal vector $\nn_{\textbf{fm}}\neq \nn$ (except for a pipe with uniform section)}
  \label{FigureSupplementaire}
 \end{center}
\end{figure}

Let $(O,\ii,\ji,\ki)$ be a convenient Cartesian reference frame, for instance the canonical basis of  $\R^3$.
The Euler equations in Cartesian coordinate are :
\begin{equation}\label{E3DICartesian}
\left\{
\begin{array}{rcl} \dive(\rho_0\,\UU) &=&0\\
\partial_t (\rho_0\,\UU) + \rho_0\,\UU\cdot\nabla \UU + \nabla P &=& \rho_0\FF
\end{array}\right.
\end{equation}
where $\UU(t,x,y,z)$ is the velocity field of components $(u,v,w)$,
$P=p(t,x,y,z)I_3$ is the isotropic pressure tensor, $\rho_0$ is the density of the water at atmospheric pressure   $p_0$ and  $\FF$ is the exterior strength of gravity given by:
$\FF=-g\vectrois{-\sin\theta(x)}{0}{\cos\theta(x)}$ where  $\theta(x)$ is the angle  $(\ii,\TT)$ in the $(\ii,\ki)$-plane (c.f. \fig\ref{FigureOxOz} or \fig\ref{FigureSupplementaire}) with $\TT$ the tangent vector (defined below) and $g$
is the gravity constant.

\noindent We introduce a characteristic function, $\phi$, in order to define the fluid area (as in \cite{GP01,M07}) :
\begin{equation}\label{DefinitionPhi}
\phi = \left\{\begin{array}{lll}
                1 & \textrm{if} & z\in \Omega(t,x),\\
                0 &  \textrm{otherwise}&
                \end{array}\right.
\end{equation}
where $\Omega(t,x)$ is the wet section  \eqref{DefinitionOmegatx}).
Using the divergence free equation, we obviously find the equation on $\phi$:
\begin{equation}\label{PhiFluidRegionAdvected}
\partial_t (\rho_0\,\phi) + \mbox{div}(\rho_0 \phi \UU) = 0.
\end{equation}
\begin{rque}
This equation  allows to get the kinematic free surface condition:
\begin{center}
{\rm any free surface particle is advected by the fluid velocity.}
\end{center}
\end{rque}
\noindent On the wet boundary ($\textbf{fm}$), we assume a no-leak condition:
$\UU\cdot\nn_{\textbf{fm}} = 0$
where $\nn_{\textbf{fm}}$ is the outward unit normal vector to the wet boundary  (as displayed on
\fig\ref{FigureSupplementaire}).
We also assume that the pressure at the free surface level is equal to the atmospheric pressure (which is assumed to be zero in the rest of the paper for the sake of simplicity).

\noindent We define the domain  $\Omega_F(t)$ of the flow at time $t$ as the union of sections,
$\Omega(t,x)$, assumed to be simply connected compact sets, orthogonal  to some  plane curve $\mathcal{C}$.
We define the parametric representation of this curve by $x\rightarrow (x,0,b(x))$  in the plane  $(O,\ii,\ji,\ki)$ where $\ki$ is the vertical axis,  $b(x)$ is the elevation of the point $\omega (x,0,b(x))$ in the $(O,\ii,\ji)$-plane (c.f. \fig\ref{FigureOxOz}).

\begin{figure}[H]
 \begin{center}
 \includegraphics[scale = 0.4]{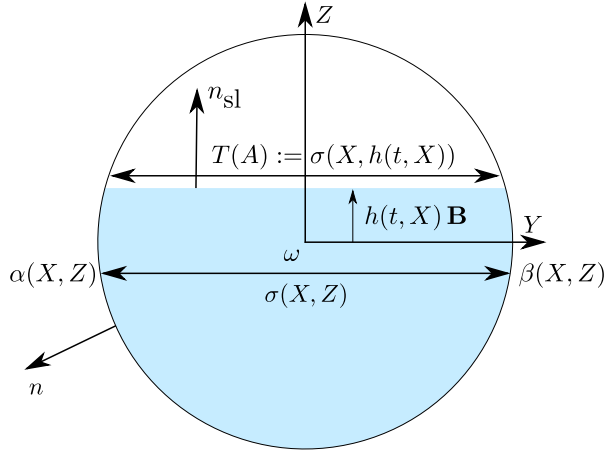}
  \caption{Transversale section $\Omega(t,X)$ at the point  $\omega$ for a free surface flow}
  \label{FigureOyOz}
 \end{center}
\end{figure}

\noindent Setting
\begin{equation}\label{VariableCurviligne1}
\dsp X = \int_{x_0}^x \sqrt{\dsp
1+\left(\frac{d\,\, b(\xi)}{dx}\right)^2} d\xi
\end{equation}
the curvilinear variable where  $x_0$ is a given abscissa, $Y = y$,
the variable ``width" and $Z$ the  $\BB$-coordinate
(i.e. the elevation of a fluid particle $M$ along the $\BB$ vector as defined below),
we define the local reference of origin  $\omega(x,0,b(x))$ and by the basis $(\TT,\NN,\BB)$ where $\TT$ is the unit tangent vector,  $\NN$ the unit normal vector and $\BB$ the unit bi-normal vector attached to the plane curve $\mathcal{C}$ at the point
$\omega(x,0,b(x))$ (see \fig\ref{FigureOxOz} and \fig \ref{FigureOyOz} for notations). In the $(O,\ii,\ki)$-plane, the vector $\BB$ is normal to the curve $\mathcal{C}$.

\noindent With these notations, for every point $\omega \in \mathcal{C}$, the wet section $\Omega(t,X)$
can be defined by the following set:
\begin{equation}\label{DefinitionOmegatx}
\Omega(t,X) = \left\{(Y,Z)\in\R^2; Z \in [-R(X),-R(X)+H(t,X)],\, Y \in [\alpha(X,Z),\beta(X,Z)]\right\}
\end{equation}
\noindent where $R(X)$ is the radius of the pipe section $S(X)=\pi R^2(X)$ and $H(t,X)$ is the physical water height.
We note $\alpha(X,Z)$ (respectively $\beta(X,Z))$ the left  (respectively  right) boundary point at elevation
$Z$, for  $-R(X)\leqslant Z\leqslant R(X)$ (as displayed on \fig\ref{FigureOyOz}).
We also assume that the support of the functions $\alpha(\cdot,z)$ and $\beta(\cdot,z)$ are compact in $[-R(X),R(X)]$.
Finally, we note the $Z$-coordinate of the water height by $h(t,X) = -R(X)+H(t,X)$.

In what follows, we will assume that the following condition holds:
\begin{description}
\item[$(H)$] Let $\mathcal{R}(x)$ bet the algebraic curvature radius at the point $\omega(x,0,b(x))$.
Then, for every  $x\in\mathcal{C}$, we have:
$$\abs{\mathcal{R}(x)}> R(x).$$
\end{description}
\begin{rque}
This geometric condition ensure that the application $\mathcal{T} : (x,y,z)\rightarrow (X,Y,Z)$  is a $C^1$-diffeomorphism. In other words, it simply means that for a given fluid particle,  there exists a unique point $\omega\in\mathcal{C}$ as displayed on \fig\ref{FigureHypotheseH}.
 \begin{figure}[H]
  \begin{center}
  \includegraphics[scale = 0.45]{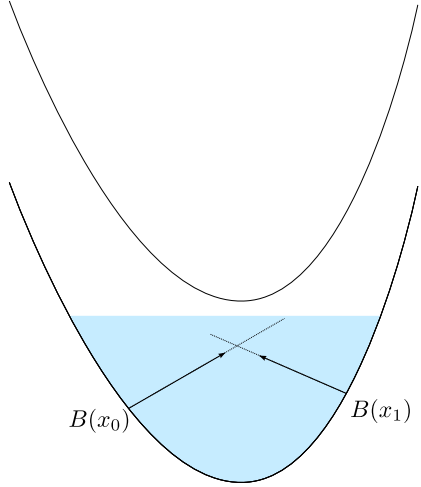}
   \caption{Forbidden case by assumption $(H)$ :}
   \label{FigureHypotheseH}
  \end{center}
 \end{figure}
\end{rque}
\subsection{Incompressible Euler equations in curvilinear coordinates}\label{SubsectionEulerCurvilinear}
Following Bouchut \emph{et al.} \cite{BFML07,BMCPV03}, we write the previous system (\ref{E3DICartesian})
in the local frame of origin $\omega(x,0,b(x))$ and of basis $(\TT,\NN,\BB)$ by the following change of variables
 $\mathcal{T} : (x,y,z)\rightarrow (X,Y,Z)$ using the divergence chain rule that we recall here:
\begin{lem}\label{lemma}

Let $(X,Y,Z) \mapsto \mathcal{T}(X,Y,Z)=(x,y,z)$ be a $C^1$-diffeomorphism  and $\mathcal{A}^{-1} =
\nabla_{(X,Y,Z)} \mathcal T$ be the Jacobian matrix with determinant $J$.

\noindent Then, for every vector field $\Phi$, we have:
$$J \dive_{(x,y,z)} \Phi = \dive_{(X,Y,Z)}(J\mathcal{A}\Phi)\,,$$
and, for every scalar function  $f$ :
$$\nabla_{(x,y,z)}f = \mathcal{A}^t \nabla_{(X,Y,Z)}f,$$ where $\mathcal{A}^t$ is the transposed matrix of $\mathcal{A}$.
\end{lem}
\noindent Let $(U,V,W)^t$ be the components of the vector field in variables $(X,Y,Z)$, $$(U,V,W)^t = \Uptheta
(u,v,w)^t$$
 where $\Uptheta$ is the rotation matrix generated  around the axis $\ji$:
$$\Uptheta=\left(
\begin{array}{ccc}
 \cos\theta & 0 & \sin\theta \\
 0 & 1 & 0 \\
 -\sin\theta& 0 & \cos\theta \\
\end{array} \right)\,.$$
\subsubsection{Transformation of the divergence equation}\label{SubsubsectionEquationEulerIncompressibleDivergenceCurvilineaire}
A given point $M$ of coordinates $(x,y,z)$ such that:
\begin{equation}\label{DefinitionMxyz}
M(x,y,z)=\Big(x- Z \sin\theta(x)\,,\,y\,,\, x+ Z \cos\theta(x)\Big)
\end{equation}
(in the $(O,\ii,\ji,\ki)$-basis) has $(X,Y,Z)$-coordinate in the local frame generated by the basis $(\TT,\NN,\BB)$ from origin $\omega$  and the matrix $\mathcal{A}^{-1}$ (appearing in Lemma \ref{lemma}) reads as follows:
$$\begin{array}{rcl}
 \mathcal{A}^{-1}&=&\left(
\begin{array}{ccc}
\dsp  \frac{d\, x}{dX} -Z \frac{d\, \theta}{dX}  \cos\theta(X) & 0 & \sin\theta(X) \\
\dsp  0 & 1 & 0 \\
\dsp  \frac{d\, b}{dX} -Z \frac{d\,\theta}{dX}  \sin\theta(X) & 0 & \cos\theta(X) \\
\end{array}
\right)\\ \\
  &=& \left(
\begin{array}{ccc}
 J\cos\theta & 0 & \sin\theta \\
 0 & 1 & 0 \\
 J\sin\theta& 0 & \cos\theta \\
\end{array}
\right)
\end{array}
$$
where
$$\dsp \frac{d\, x}{dX}  = \frac{1}{\dsp \sqrt{1+\left(\frac{d\,b}{dx} \Big(x(X,Y,Z)\Big)\right)^2}}= \cos\theta(X),
\quad
\frac{d\,b(X)}{dX}  = \sin\theta(X),$$
and $$J(X,Y,Z):=\det(\mathcal{A}^{-1}) = 1-Z \frac{d\,\theta(X)}{dX}
$$
with $J(X,Y,Z)=J(X,Z)$.

\noindent Then, we have:
\begin{equation}\label{DefinitiondeAmoin1}
\mathcal{A}=\frac{1}{J}\left(
\begin{array}{ccc}
 \cos\theta & 0 & \sin\theta \\
 0 & 1 & 0 \\
 -J\sin\theta& 0 & J\cos\theta \\
\end{array}
\right)
\end{equation}
and using Lemma (\ref{lemma}), the free divergence equation in variables $(X,Y,Z)$ is,
$$J\dive_{x,y,z}(\UU) = \dive_{X,Y,Z}\left(\begin{array}{c} U\\JV\\JW \end{array}\right)=0,$$i.e.
\begin{equation}\label{Incomp}
\partial_X U + \partial_Y(JU)+\partial_Z(JW) = 0 .
\end{equation}
\begin{rque}
The application $(x,y,z)\rightarrow M(x,y,z)$ is a $C^1$-diffeomorphism since $J(X,Z)>0$  in view of the assumption $(H)$.
\end{rque}

\subsubsection{Transformation of the equation of  conservation of the momentum}\label{SubsubsectionEquationEulerIncompressibleMomentCurvilineaire}
Following the previous paragraph, using Lemma \ref{lemma} to the scalar convection equation $\dsp\frac{d f}{dt} $
characterized by the speed $\UU$ which is a divergence free field, we get the following identity:
$$J(\partial_t+\UU\cdot\nabla )f =
J\Big(\partial_t f+\dive(f \UU)\Big)=J \dive_{t,x,y,z} \vecdeux{J f}{J\mathcal{A}^{-1}f\UU},$$
where $\mathcal{A}^{-1}$ is the inverse matrix of $\mathcal{A}$ given by (\ref{DefinitiondeAmoin1}).
Thus, we have:
\begin{equation}\label{ET}
\begin{array}{lll}J(\partial_t+\UU.\nabla )f &=&
\partial_t(Jf)+\partial_X(fU)+\partial_Y(JfV)+\partial_Z(JfW)
\end{array}
\end{equation}

Performing a left multiplication  of the equation of conservation of the momentum \eqref{E3DICartesian} by $J \Uptheta$, where the source term is written as $\FF=-\nabla \left(\textbf{g}\cdot
M\right)$ (for a point $M$ defined as previously \eqref{DefinitionMxyz}), we get:
$$\begin{array}{lll}
0&=&J\Uptheta(\partial_t \UU + \UU\cdot\nabla \UU + \dive(P/\rho_0) +\nabla \left(\textbf{g}\cdot M\right) \\
\vspace{0.8mm}
 &=&J\Big(\partial_t (\Uptheta\UU) + (\Uptheta\UU\cdot\nabla) \UU + J\Uptheta \dive(P/\rho_0)+J\Uptheta\nabla
\left(\textbf{g}\cdot M\right)\Big)\\ \vspace{0.8mm}
 &=& \underbrace{J\left(\partial_t \vectrois{U}{V}{W} + \vectrois{(\UU\cdot\nabla u) \cos\theta+(\UU\cdot\nabla w)
\sin\theta}{\UU\cdot\nabla v}{-(\UU\cdot\nabla u)\sin\theta+(\UU\cdot\nabla w)\cos\theta}\right)}_{(a)}  \\
& &+
\underbrace{\vectrois{J\dive(\psi\ii)\cos\theta+J \dive(\psi\ki)\sin\theta}{J\dive(\psi\ji)}{
-J\dive(\psi\ii)\sin\theta+J\dive(\psi\ki)\cos\theta}}_{(b)}
\end{array}$$
where $\psi:=(p+g(b+Z\cos\theta))/\rho_0$.

\noindent Then, we proceed in two steps:
\paragraph{Computation of  $(a)$.\newline}
\noindent We have:
\begin{equation}\label{terme_a}
\begin{array}{l}
J\vectrois{\partial_t U  + (\UU\cdot\nabla u) \cos\theta+(\UU\cdot\nabla w)\sin\theta}
          {\partial_t V  + \UU\cdot\nabla v}
          {\partial_t W  + -(\UU\cdot\nabla u)\sin\theta+(\UU\cdot\nabla w)\cos\theta} \\
= J\vectrois{\partial_t U  + \UU\cdot\nabla U-W \UU\cdot\nabla \theta}
             {\partial_t V  + \UU\cdot\nabla V }
             {\partial_t W  + \UU\cdot\nabla W+U\UU\cdot\nabla \theta}.
\end{array}
\end{equation}
Applying successively the identity \eqref{ET} with $f=U,V,W$, we get:
\begin{equation}\label{terme(a)}
\partial_t \vectrois{JU}{JV}{JW} + \dive_{X,Y,Z}\left(\vectrois{U}{JV}{JW}\otimes\vectrois{U}{JV}{JW}\right)
-\ii UW\frac{d\, \theta}{dX} + \ki U^2\frac{d\, \theta}{dX}.
\end{equation}
\paragraph{Computation of  $(b)$.\newline}
\noindent Applying again Lemma \ref{lemma}, we show that the three following identities hold for every scalar function $\psi$:
\begin{equation}\label{identite1}
\left\{
\begin{array}{lll}
J \dive(\psi \ii) &=& \dive_{X,Y,Z} \vectrois{\psi\cos\theta}{0}{-J\psi\sin\theta},\\
J \dive(\psi \ji) &=&\partial_Y(J\psi),\\
J \dive(\psi \ki) &=& \dive_{X,Y,Z} \vectrois{\psi\sin\theta}{0}{J\psi\cos\theta}.
\end{array}
\right.
\end{equation}
Moreover, we have:
\begin{equation}\label{identite2}
\left\{\begin{array}{lll}
\partial_X(\psi\cos\theta)\cos\theta+\partial_X(\psi\sin\theta)\sin\theta
&=& \partial_X \psi,\\
\partial_Z(J\psi\cos\theta)\sin\theta-\partial_Z(J\psi\sin\theta)\cos\theta&=&0,
\end{array}\right.
\end{equation}
and
\begin{equation}\label{identite3}
\left\{\begin{array}{lll}
\partial_X(\psi\sin\theta)\cos\theta-\partial_X(\psi\cos\theta)\sin\theta
&=& \psi\partial_X \theta,\\
\partial_Z(J\psi\cos\theta)\cos\theta+\partial_Z(J\psi\sin\theta)\sin\theta&=&\partial_Z(\psi J).
\end{array}
\right.
\end{equation}
In view of equalities \eqref{identite1}--\eqref{identite3} applied to the quantity $\psi:=(p+g(b+Z\cos\theta))/\rho_0$,
the term $(b)$ reads as follows:
\begin{equation}\label{terme(b)}
\vectrois{\partial_X(\psi)}{\partial_Y(J\psi)}{\psi\partial_X\theta +\partial_Z(J\psi)} .
\end{equation}
Finally, gathering results \eqref{terme(a)}--\eqref{terme(b)}, the incompressible Euler equations in
variables $(X,Y,Z)$ are:
\begin{equation}\label{E3DICurvilinear}\left\{
\begin{array}{lll}
\partial_X (\rho_0\,U) + \partial_Y(J\rho_0\,V)+  \partial_Z(J\rho_0\,W)&=&0,\\
\partial_t(J\rho_0\,U) + \partial_X(\rho_0\,U^2) + \partial_Y(J\rho_0\,UV)+  \partial_Z(J\rho_0\,UW)+ \partial_X p
&=& G_1,\\
\partial_t(J\rho_0\,V) + \partial_X(\rho_0\,UV) + \partial_Y(J\rho_0\,V^2)+  \partial_Z(J\rho_0\,VW)+
J\partial_Y(p) &=& 0,\\
\partial_t(J\rho_0\,W) + \partial_X(\rho_0\,UW) + \partial_Y(J\rho_0\,VW)+  \partial_Z(J\rho_0\,W^2)+
J\partial_Z(p) &=& G_2
\end{array}\right.
\end{equation}
where
\begin{equation}\label{G1G2}
G_1 =
\rho_0\,UW \frac{d\, \theta}{dX}  -g\rho_0 J\,\sin\theta,\,\,  G_2=-\rho_0\,U^2\frac{d\, \theta}{dX}
-Jg\rho_0\,\cos\theta.
\end{equation}

\noindent The no-leak condition, with respect to the new variables, becomes:
\begin{eqnarray}
\vectrois{U}{V}{W}\cdot \nn_\textbf{fm}=0\, \label{NoLeakCond}
\end{eqnarray}
where $\nn_\textbf{fm} = \dsp\frac{1}{\cos\theta(X)}\vectrois{-\sin\theta(X)}{0}{\cos\theta(X)}$. \\
The condition at the free surface, in the new variables, reads:
\begin{eqnarray}
p\Big(t,X,Y,Z=h(t,X)\Big)=0\, .\label{FSCondition}
\end{eqnarray}

\subsubsection{Formal asymptotic analysis}\label{SectionAnalysAsymptotiqueFormelleSL}
Taking advantage of the ratio of the domain, we perform  a formal asymptotic analysis of the
equations  \eqref{E3DICurvilinear} with respect to a small parameter $\varepsilon$. Especially, we are interested on the approximation at main order. In that case, we will get $J \equiv 1$.

\noindent To this end, let $\; \dsp \frac{1}{\varepsilon} = \frac{L}{H}= \frac{L}{l},
$ be the ratio aspect of the pipe, assumed   very large.  $H$, $L$ and $l$, is the characteristic height, length  and width (to simplify, $l=H$ as the pipe, here, is assumed with circular cross section).
In the same way, denoting by ($\overline{V}$,$\overline{W}$) the characteristic speed following the normal and bi-normal direction, $\overline{U}$ the characteristic speed following the main pipe axis, we also assume that:
$$\epsilon = \frac{\overline{V}}{\overline{U}}= \frac{\overline{W}}{\overline{U}}.$$
Let $T$ and $P$ be  the characteristic time and pressure such that $$\overline{U}=\frac{L}{T},\quad P = \rho_0
\overline{U}^2 .$$
We set the following non-dimensioned variables:
$$\widetilde{U} = \frac{U}{\overline U},\, \widetilde V = \varepsilon \frac{V}{ \overline U},\,
\widetilde W = \varepsilon \frac{W}{\overline U},\, $$ $$\widetilde X = \frac{X}{L},\, \widetilde Y =
\frac{Y}{H},\,
\widetilde Z = \frac{Z}{H},\,
\widetilde p = \frac{p}{P},\,\widetilde\theta = \theta.$$
\noindent Under these assumptions, the rescaled Jacobian is:
$$\widetilde{J} (\widetilde X,\widetilde Y,\widetilde Z)
= 1- \varepsilon \widetilde Z\dsp \frac{d\, \widetilde{\theta}}{d\widetilde{X}}.$$

\noindent Then, the non-dimensioned system  \eqref{E3DICurvilinear} is reduced to:
\begin{equation}\label{E3DICurvilinearRescaled}\left\{
\begin{array}{rcl}
\partial_{\widetilde X} \widetilde U + \partial_{\widetilde Y}(\widetilde{J} \widetilde V)+  \partial_{\widetilde
Z}
(\widetilde{J} \widetilde W)&=&0\\
\partial_{\widetilde t}(\widetilde{J}\, \widetilde U) + \partial_{\widetilde X}({\widetilde U}^2) +
\partial_{\widetilde Y}(\widetilde{J}\, \widetilde U\, \widetilde V)+
\partial_{\widetilde Z}(\widetilde{J}\, \widetilde U\, \widetilde W)+ \partial_{\widetilde X} \widetilde p &=&
G_{1},\\
\varepsilon^2\left(\partial_{\widetilde t}(\widetilde{J}\, \widetilde V) +
\partial_{\widetilde X}(\widetilde U\, \widetilde V) + \partial_{\widetilde Y}
(\widetilde{J}\, {\widetilde V}^2)+  \partial_{\widetilde Z}(\widetilde{J}\, \widetilde V\, \widetilde W)\right)+
\partial_{\widetilde Y}(\widetilde{J}\, \widetilde p) &=& 0,\\
\varepsilon^2\left(\partial_{\widetilde t}(\widetilde{J}\, \widetilde W) +
\partial_{\widetilde X}(\widetilde U\, \widetilde W) + \partial_{\widetilde
Y}(\widetilde{J}\, \widetilde V\, \widetilde W) +  \partial_{\widetilde Z}(\widetilde{J}\,
{\widetilde W}^2)\right)  + \widetilde{J}\partial_{\widetilde Z}(\widetilde p) &=&G_{2}
\end{array}\right.
\end{equation}
\noindent where
$$G_{1} =  \varepsilon \widetilde U \widetilde W \frac{d\, \widetilde{\theta}}{d\widetilde{X}}
-\dsp\frac{\sin\widetilde\theta(\widetilde X)}{{F_{r,L}}^2}-
  \frac{\widetilde Z \dsp }{{F_{r,H}}^2} \frac{d\,}{d\widetilde{X}}\cos\widetilde{\theta}(\widetilde X), $$
$$G_{2} =
-\varepsilon{\widetilde U}^2 \frac{d\, \widetilde{\theta}}{d\widetilde{X}}-
\dsp\frac{\cos\widetilde\theta(\widetilde X)}{{F_{r,H}}^2}
+
\dsp\varepsilon\frac{d\, \widetilde{\theta}}{d\widetilde{X}}\frac{\widetilde Z\widetilde{J}\cos\widetilde\theta(\widetilde X)}{{F_{r,H}}^2},$$
 $F_{r,\chi} = \dsp\frac{\overline{U}}{\sqrt{g \chi}}$ is the Froude number following the axis $\TT$, $\BB$ or
$\NN$ with $\chi=L$  or $\chi=H$.

\noindent Formally, taking  $\varepsilon =0 $ in the previous equation,  we get:
\begin{eqnarray}
\partial_{\widetilde X} \widetilde U + \partial_{\widetilde Y} \widetilde V+  \partial_{\widetilde Z} \widetilde
W&=&0, \label{E3DICurvilinearRescaledEpsilonVanishes1} \\
\partial_{\widetilde t} \widetilde{U} + \partial_{\widetilde X} ( {\widetilde U}^2) + \partial_{\widetilde Y} (
\widetilde U  \widetilde V)+ \partial_{\widetilde Z} ( \widetilde U  \widetilde W)+ \partial_{\widetilde X}
\widetilde p &=&
  -\dsp\frac{\sin\widetilde\theta(\widetilde X)}{{F_{r,L}}^2}  - \dsp  \frac{\widetilde Z \dsp}{{F_{r,H}}^2}
\frac{d\,}{d\widetilde X}
\cos\widetilde\theta(\widetilde X)
\label{E3DICurvilinearRescaledEpsilonVanishes2},\\
\partial_{\widetilde Y} \widetilde p &=& 0 \label{E3DICurvilinearRescaledEpsilonVanishes3}, \\
\partial_{\widetilde Z} \widetilde p &=& -\dsp\frac{\cos\widetilde\theta(\widetilde X)}{{F_{r,H}}^2}
\label{E3DICurvilinearRescaledEpsilonVanishes4}.
\end{eqnarray}

\noindent\textcolor{red}{\danger} Henceforth, we note $(x,y,z)$ the dimensioned variables $(X,Y,Z)$ and $(u,v,w)$ dimensioned speed $(U,V,W)$.
In particular, we set:
$$x=L \widetilde X,\,y=H \widetilde Y,\,z= H \widetilde Z$$ et
$$u=\overline{U} \widetilde U,\,v=\epsilon \overline{U} \widetilde V ,\,v=\epsilon \overline{U} \widetilde W \textrm{ et } p = \rho_0 \overline{U}^2\, \widetilde p.$$
Then, multiplying the equation (\ref{E3DICurvilinearRescaledEpsilonVanishes1}) by $ \dsp \rho_0\frac{\overline{U}}{L}$, (\ref{E3DICurvilinearRescaledEpsilonVanishes2}) by $ \dsp \rho_0\frac{\overline{U}}{T}$,
 (\ref{E3DICurvilinearRescaledEpsilonVanishes3}) by $\dsp\rho_0\frac{\overline{U}^2}{H}$,
 (\ref{E3DICurvilinearRescaledEpsilonVanishes4}) by $\dsp\rho_0\frac{\overline{U}^2}{H}$,
we obtain the  hydrostatic approximation of the Euler equations (\ref{E3DICurvilinear}):
\begin{eqnarray}
\partial_x (\rho_0 u) + \partial_y (\rho_0 v)+  \partial_z (\rho_0 w)&=&0,
\label{E3DICurvilinearRescaledDimensionnelles1} \\
\hspace{-7mm}\partial_t (\rho_0 u) + \partial_x ( \rho_0 u^2) + \partial_y (\rho_0 u v)+ \partial_z ( \rho_0  u w)+ \partial_x p
&=&   -\dsp g\rho_0 \sin\theta(x)  -  g \rho_0 z \frac{d\,}{dx}\cos\theta(x)
\label{E3DICurvilinearRescaledDimensionnelles2},\\
\partial_y p &=& 0 \label{E3DICurvilinearRescaledDimensionnelles3}, \\
\partial_z p &=& -g \cos\theta(x) \label{E3DICurvilinearRescaledDimensionnelles4}.
\end{eqnarray}

\subsubsection{Vertical averaging of the hydrostatic approximation of Euler equations}
\label{SubsubsectionMoyennisationEulerIncompressibleCurvilineaire}
Let  $A(t,x)$ and $Q(t,x)$ be the conservative variables of wet area and discharge defined by the follwoing relations:
\begin{equation}\label{DefinitionA}
A(t,x) = \int_{\Omega(t,x)}\, dydz
\end{equation}
and
\begin{equation}\label{DefinitionQ}
Q(t,x) = A(t,x) \overline{u}(t,x)
\end{equation}
where \begin{equation}\label{DefinitionuEulerIncomp}
\overline{u}(t,x) = \frac{1}{A(t,x)}\int_{\Omega(t,x)} u(t,x,y,z)\, dydz
\end{equation}
is the mean speed of the fluid over the section $\Omega(t,x)$.

\paragraph{Kinematic boundary condition and the equation of the conservation of the mass.\newline}
Let $\dsp \textbf{V}$ be the vector field $\vecdeux{v}{w}$. Integrating the equation of conservation of the mass \eqref{PhiFluidRegionAdvected} on the set:
$$\overline{\Omega}(x) = \{(y,z);\, \alpha(x,z)\leqslant y \leqslant \beta(x,z),\,  -R(x) \leqslant y
\leqslant \infty  \},$$ we get the follwing equation:
\begin{equation}\label{CLSL1}
\int_{\overline{\Omega}(x)} \partial_t (\rho_0 \phi) + \partial_x (\rho_0 \phi u)+ \dive_{y,z} (\rho_0 \phi
\textbf{V}) \,dy dz = \rho_0\left(\partial_t A + \partial_x Q + \int_{\partial\overline{\Omega}_{\textrm{fm}}(x)}
\left(u \partial_x M-\textbf{V}\right) \cdot \nn\,ds   \right)
\end{equation}
where $A$ and $Q$ are given by \eqref{DefinitionA} and \eqref{DefinitionQ}.

\noindent According to the  definition \eqref{DefinitionPhi} of  $\phi$, the boundary
$\overline{\Omega}_{\textrm{fm}}$ coincides with $\gamma_{\textrm{fm}}$.
Using, the no-leak condition \eqref{NoLeakCond}, Equation \eqref{CLSL1} is  equivalent to
\begin{equation}\label{CLSL2}
\partial_t (\rho_0 A) + \partial_x(\rho_0 Q) = 0
\end{equation}

\noindent Now, if integrate the equation \eqref{PhiFluidRegionAdvected} on $\Omega(t,x)$, we get:
\begin{equation}\label{CLSL3}
\begin{array}{l}
\rho_0\left(\dsp \int_{-R(x)}^{h(t,x)} \partial_t \int_{\alpha(x,z)}^{\beta(x,z)}\,dy dz + \partial_x Q +
\int_{\partial\Omega(t,x)} \left(\textbf{V}-u\partial_x M\right)\cdot\nn \,ds\right)=0
\end{array}
\end{equation}
where
$$
\begin{array}{l}
\dsp \int_{-R(x)}^{h(t,x)} \partial_t \int_{\alpha(x,z)}^{\beta(x,z)}\,dy dz = \partial_t A  - \sigma(x,h(t,x))
\partial_t h
\end{array}
$$
with  $\sigma(x,h(t,x))$ the width at the free surface elevation as displayed on \fig\ref{FigureOyOz}.

\noindent In view of the no-leak condition \eqref{NoLeakCond}, the integral on the wet boundary is zero, i.e. :
$$\int_{\gamma_{\textrm{fm}}(t,x)} \left(\textbf{V}-u\partial_x M\right)\cdot\nn_{\textrm{fm}} \,ds=0 .$$
Then, we deduce:
\begin{equation}\label{CLSL4}
\begin{array}{l}
\dsp \partial(\rho_0 A) + \partial_x(\rho_0 Q) + \rho_0\int_{\gamma_{\textrm{sl}}(t,x)} \left(\partial_t M
+ u\partial_x
M - \textbf{V}\right)\cdot\nn_{\textrm{sl}} \,ds = 0.
\end{array}
\end{equation}
By comparing equations (\ref{CLSL2}) and (\ref{CLSL4}), we finally get the kinematic condition at the free surface:
\begin{equation}\label{ConditionALaSurfaceLibre}
\int_{\gamma_{\textrm{sl}}(t,x)} \left(\partial_t M + u\partial_x
M - \textbf{V}\right)\cdot\nn_{\textrm{sl}} \,ds =0 .
\end{equation}
and we  deduced from \eqref{CLSL4} the following equation of the conservation of the mass:
\begin{equation}\label{SLMasse}
\partial_t (\rho_0 A) + \partial_x (\rho_0  Q) =0 .
\end{equation}

\paragraph{Equation of the conservation of the momentum.\newline}
In order to get the equation of the conservation of the momentum of the free surface model, we integrate each terms of
\eqref{E3DICurvilinearRescaledDimensionnelles2} along sections $\Omega(t,x)$ as follows:
$$\dsp\int_{\Omega} \underbrace{\partial_t (\rho_0 u)}_{a_1} +
\underbrace{\partial_x(\rho_0  u^2)}_{a_2} +
\underbrace{\dive_{y,z} \left(\rho_0 u \textbf{V}\right)}_{a_3}+
\underbrace{\partial_x p}_{a_4} \,dy dz= \int_{\Omega}
\underbrace{-\rho_0 g z  \frac{d\,}{dx}\cos\theta}_{a_5}-\underbrace{\rho_0  g\sin\theta}_{a6} \,dy dz$$
where $\dsp \textbf{V} = \vecdeux{v}{w}$. \\
Assuming that
$$\overline{u\,v} \approx \overline{u} \,\overline{v},\, \overline{u^2} \approx \overline{u}^2,$$
we successively get:
\paragraph{Computation of the term $\bm{ \dsp \int_{\Omega(t,x)} a_1\,dy dz}$.\newline}
The pipe being non-deformable, only the integral at the free surface is relevant:
$$\dsp \int_{\gamma_{\textrm{fm}}(t,x)} \rho_0  u \,\partial_t M \cdot \nn_{\textrm{sl}} \,ds = 0.$$
So, we get:
$$\begin{array}{lll}
\dsp\int_{\Omega(t,x)}{\partial_t(\rho_0 u)} \,dy dz &=& \dsp\partial_t \int_{\Omega(t,x)} \rho_0 u \, dydz -
\int_{\gamma_{\textrm{sl}}(t,x)} \rho_0  u \,\partial_t M \cdot \nn_{\textrm{sl}} \,ds .
\end{array}
$$
\paragraph{Computation of the term $\dsp \bm{\int_{\Omega(t,x)} a_2\,dy dz}$.\newline}
$$\begin{array}{lll}
\dsp\int_{\Omega(t,x)}{\partial_x(\rho_0  u^2)} \,dy dz &=& \dsp\partial_x \int_{\Omega(t,x)} \rho_0  u^2 \,dydz -
\int_{\gamma_{\textrm{sl}}(t,x)}\rho_0 u^2 \partial_x M \cdot \nn_{\textrm{sl}} \,ds\\
& & \dsp -\int_{\gamma_{\textrm{fm}}(t,x)}\rho_0 u^2 \partial_x M \cdot \nn_{\textrm{fm}} \,ds.
\end{array}
$$
\paragraph{Computation of the term $\dsp \bm{\int_{\Omega(t,x)} a_3\,dy dz}$.\newline}
$$\begin{array}{lll}
\dsp\int_{\Omega(t,x)}\dive_{y,z} \left(\rho_0 u \textbf{V}\right)\, dy dz
&=& \dsp \int_{\gamma_{\textrm{sl}}(t,X)} \rho_0 u \textbf{V}\cdot\nn_{\textrm{sl}}\, ds\\
& & \dsp +\int_{\gamma_{\textrm{fm}}(t,X)} \rho_0 u \textbf{V}\cdot\nn_{\textrm{fm}}\, ds.
\end{array}$$
Summing  the result of the previous step $a_1+a_2+a_3$, we get:
\begin{equation}\label{Etape1-2-3}
\dsp\int_{\Omega(t,x)} a_1+a_2+a_3 \,dy dz = \dsp\partial_t(\rho_0 Q) + \partial_x\left(\rho_0\frac{Q^2}{A}\right)
\end{equation} where $A$ and $Q$ are given by (\ref{DefinitionA}) and (\ref{DefinitionQ}).
\paragraph{Computation of the term $\dsp \bm{\int_{\Omega(t,x)} a_4\,dy dz}$.\newline}
Let us first note that the pressure is   hydrostatic:
\begin{equation}\label{PressionHydrostatique}
p(t,x,z) = \rho_0 g(h(t,x)-z)\cos\theta(x)
\end{equation}
since from the equation \eqref{E3DICurvilinearRescaledDimensionnelles3}, the pressure does not depend on the variable $y$. Equation \eqref{PressionHydrostatique} follows immediately by integrating the equation (\ref{E3DICurvilinearRescaledDimensionnelles4}) from $z$ to $h(t,x)$.

\noindent For $\psi=p$, $p$ given by the relation  \eqref{PressionHydrostatique},  $(t,x)$ fixed, we have:
$$\begin{array}{lll}
\dsp \int_{\Omega(t,x)} \partial_x \psi\, dy dz &=&
\dsp\int_{-R(x)}^{h(t,x)}\int_{\alpha(x,z)}^{\beta(x,z)} \partial_x \psi \, dy dz\\
&=& \dsp\int_{-R(x)}^{h(t,x)}\partial_x \int_{\alpha(x,z)}^{\beta(x,z)} \psi \, dy dz\\
& &\dsp - \left(
\int_{-R(x)}^{h(t,x)}
       \partial_x \beta(x,z)\,\psi_{|y=\beta(x,z)}-\partial_x \alpha(x,z) \psi_{|y=\alpha(x,z)} \,dz
  \right)\\ &=&
\dsp \partial_x \int_{\Omega(t,x)} \psi\, dy dz  \\
& & \dsp - \left(
\int_{-R(x)}^{h(t,x)}
       \partial_x \beta(x,z)\,\psi_{|y=\beta(x,z)}-\partial_x \alpha(x,z) \psi_{|y=\alpha(x,z)} \,dz
\right) \\
& & \dsp -  \partial_x h(t,x) \int_{\alpha_{|z=h(t,x)}}^{\beta_{|z=h(t,x)}} \psi_{|z=h(t,x)} \, dy \\
& & \dsp - \partial_x R(x)   \int_{\alpha_{|z=h(t,x)}}^{\beta_{|z=h(t,x)}} \psi_{|z=-R(x)}  \, dy .
\end{array}
$$
Finally,  we have:
\begin{equation}\label{Etape4bis}
\begin{array}{lll}
\dsp\int_{\Omega(t,x)} \partial_x  p   \,dy dz &=&  \dsp \partial_x
(\rho_0 g I_1(x,A(t,x))\cos\theta(x))-g \rho_0  I_2(x,A)\cos\theta(x)\\
& & -\dsp \rho_0 g \big(h(t,x)+R(x)\big)\cos\theta(x) \sigma\big(x,-R(x)\big)\frac{d\,\, R(x)}{dx}
\end{array}
\end{equation}
where $I_1$ is the  hydrostatic pressure:
\begin{equation}\label{DefinitionI1}
I_1(x,A) = \int_{-R(x)}^{h(A)}(h(A)-z) \sigma(x,z) \,dz.
\end{equation}
When the sections of the pipe are rectangular and uniform, we have $I_1(x,A):=I_1(A)$ and
$\sigma(x,z)=\sigma=cte$. Moreover, we have $A=(h+R) \sigma= H  \sigma$ and the pressure reads
$$\dsp  \frac{g I_1(A)}{ \sigma}= \frac{g I_1(A)}{ \sigma} = g \dsp\frac{H^2}{2}$$ as for the usual formulation of the mono-dimensional Saint-Venant equations.

\begin{figure}[!ht]
 \begin{center}
 \includegraphics[scale = 0.4]{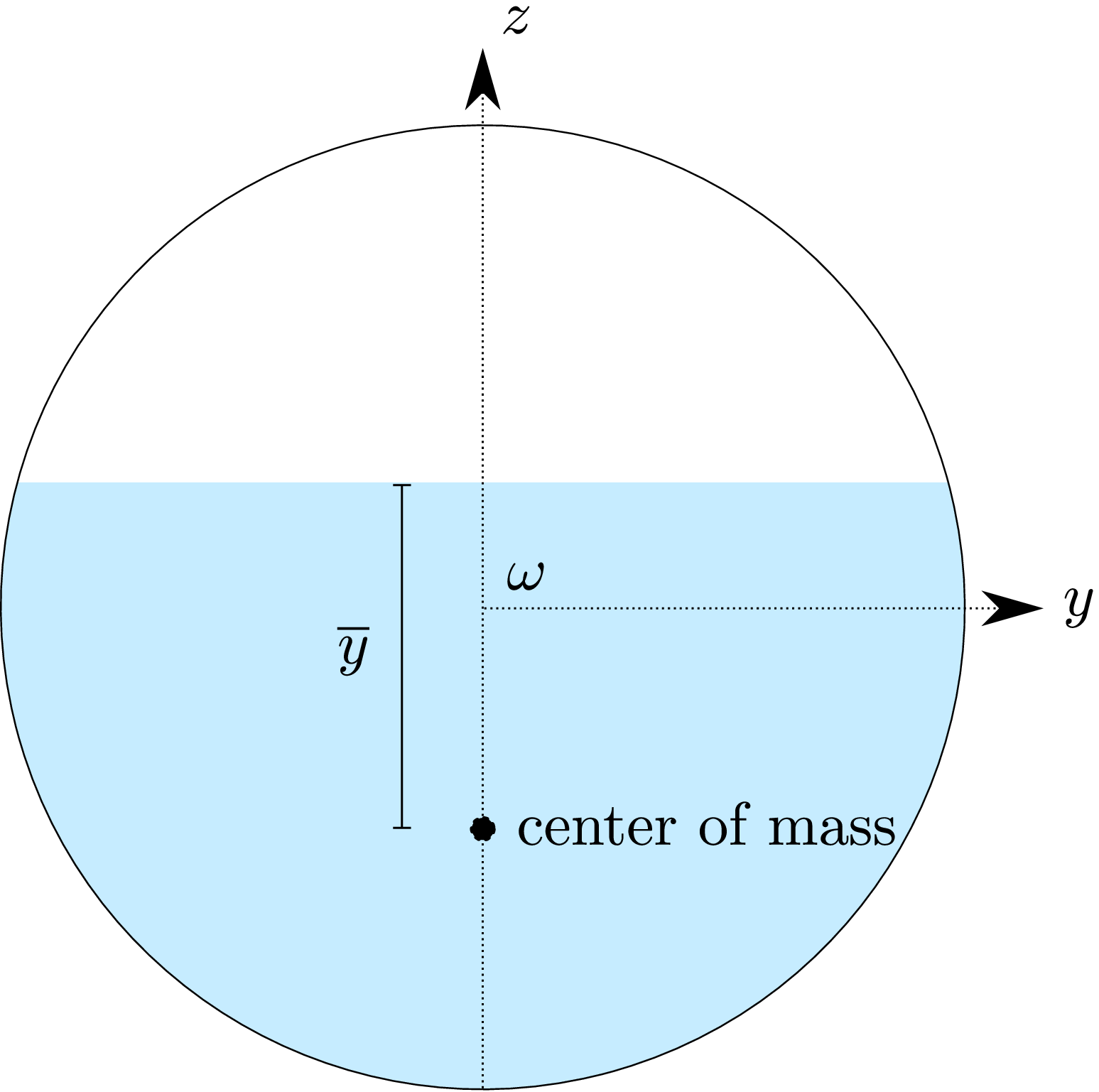}
  \caption{$\bar y$}
  \label{ybar}
 \end{center}
\end{figure}

We can also regard $I_1/A =   \overline{y}$ as the distance separating the free surface to the center of the mass of the wet section (see \fig\ref{ybar}).

\noindent The term $I_2$ is the pressure source term:
\begin{equation}\label{DefinitionI2}
I_2(x,A) = \int_{-R(x)}^{h(A)} (h(A)-z) \partial_x \sigma(x,z) \,dz.
\end{equation}
It takes into account of the section variation via the term $\partial_x \sigma(x,\cdot)$.

\noindent The term $\dsp \rho_0 g \big(h(t,x)+R(x)\big)\cos\theta(x)
\sigma\big(x,-R(x)\big)\frac{d\, R(x)}{dx}$ is also a term which takes into account the variations of the section.
The contribution of this term is non zero when:
\begin{eqnarray}
\sigma(x,z=-R(x))&\neq& 0,\textrm{ and } \label{ConditionTypeDeConduite1}\\
\partial_x R(x)  &\neq& 0.\label{ConditionTypeDeConduite2}
\end{eqnarray}
As we only deal with pipe with circular section, therefore, the result of the computation is simply $\dsp \int_{\Omega(t,x)} a_4\,dy dz$
\begin{equation}\label{Etape4}
\begin{array}{lll}
\dsp\int_{\Omega(t,x)} \partial_x  p   \,dy dz &=&  \dsp \partial_x
(\rho_0 g I_1(x,A(t,x))\cos\theta(x))-g \rho_0  I_2(x,A)\cos\theta(x)
\end{array}
\end{equation}
in the rest of the paper.
\paragraph{Computation of the term $\dsp \bm{\int_{\Omega(t,x)} a_5\,dy dz}$.\newline}
We have:
\begin{equation}\label{Etape5}
\int_{\Omega(t,x)}  \rho_0 g z  \frac{d\,}{dx}\cos\theta  \,dy dz = \rho_0 g A \overline{z}
\frac{d\,}{dx}\cos\theta
\end{equation}
where $ \overline{z}$ is the $z$-coordinate of the center of the mass.
As  $\dsp \frac{I_1(x,A(t,x)}{A(t,x)}:=\overline{y}$ (see step  ``Computation of the term $a_3$."),
the quantity $\overline{z}$ is related to  $I_1$ by the formula:
\begin{equation}\label{DefinitionDeZbar}
\overline{z} =  h(t,x)-\frac{I_1(x,A(t,x)}{A(t,x)}.
\end{equation}
\paragraph{Computation of the term $\dsp \bm{\int_{\Omega(t,x)} a_6\,dy dz}$.\newline} We have:
\begin{equation}\label{Etape6}
\int_{\Omega(t,x)}   \rho_0  g\sin\theta \,dy dz = \rho_0 g A \sin\theta.
\end{equation}

Then, gathering results (\ref{Etape1-2-3})--(\ref{Etape6}), we get the equation of the conservation of the momentum.
Finally, the new shallow water equations for free surface flows in closed water pipe with variable slope and section are:
\begin{equation}\label{ModelSL}
\hspace{-0.5mm}\left\{\begin{array}{lll}
\dsp \partial_t(\rho_0 A) +\partial_x(\rho_0 Q) &=& 0\\
\dsp \partial_t(\rho_0 Q) +\partial_x\left(\rho_0 \frac{Q^2}{A}+g \rho_0 I_1 \cos\theta\right)&=& \dsp
 -g \rho_0 A \sin\theta + g \rho_0 I_2 \cos\theta -g \rho_0 A \overline{z} \frac{d\,}{dx}\cos\theta
\end{array}\right.
\end{equation}
This model is called \FS-model.

In System (\ref{ModelSL}), we may add a friction term $-\rho_0 g S_f\, \TT$ to take into account dissipation of energy. We have chosen this term $S_f$ as the one  given by the Manning-Strickler law (see e.g. \cite{SWB98}): $$S_f(A,U)=K(A)U|U|\,.$$ The term $K(A)$ is defined by: $\dsp K(A) = \frac{1}{K_s^2 R_h(A)^{4/3}}$, $K_s>0$ is the Strickler coefficient of  roughness depending on the material,  $R_h(A)= A/P_m$ is the hydraulic radius and $P_m$ is the perimeter of the wet surface area (length of the part of the channel's section in contact with the water).

\section{Formal derivation of the pressurized model}\label{SectionDerivationEquationsP}
When the section is completely filled, we have to define a strategy to derive a suitable pressurized model in order to
\begin{itemize}
\item take into the compressibility of the water,
\item modelise the water hammer (issuing form the overpressure and depression waves)
\end{itemize}
keeping in mind that we want to construct a mixed model which allows 
\begin{itemize}
\item deal with free surface flows,
\item deal with pressurized flows and
\item to cope the transition between a free surface state and a pressurized (i.e. compressible) state transition phenomenon.
\end{itemize}
There exists a large literature on this topic, for instance
\begin{itemize}
\item the Preissmann slot artefact (see, for instance, \cite{Preissmann71})  but this technique has the drawback to do not take into account the sub-atmospheric flows,
\item the Allievi equations (see, for instance,\cite{B00}) but this equation are not well suited for a coupling with the derived \FS-model.
\end{itemize}

Then, as a starting point, we consider the $3$D isentropic compressible Euler equations:
\begin{equation}\label{Euler3D_mass_conservation}
\partial_t \rho + \dive{(\rho \UU)} = 0,
\end{equation}
\begin{equation}\label{Euler3D_momentum_conservation}
\partial_t (\rho\UU) + \dive{(\rho \UU\otimes \UU)} + \nabla p(\rho) = \rho\FF,
\end{equation}
where $\UU(t,x,y,z)$ is the fluid velocity of components $(u,v,w)$ and $\rho(t,x,y,z)$ is the volumetric mass of the fluid.
The gravity  source term is
$\FF=-g\vectrois{-\sin\theta(x)}{0}{\cos\theta(x)}$ where $\theta(x)$ is the angle $(\ii,\TT)$ (see
\fig\ref{FigureOxOz} or \fig\ref{FigureSupplementaire}).
As defined previously, $\TT$ is the tangential vector at the point $\omega\in\mathcal{C}$ (see Section \ref{SectionFormalDerivationFSModel} for notations) where the ``pressurized" plane curve is defined below.

\noindent The system is closed by the linearised pressure law (see \cite{SWB98,WS78}):
\begin{equation}\label{PressionLinearisee}
p = p_a + \dsp{\frac{\rho-\rho_0}{\beta_0 \rho_0}}
\end{equation}
which have the advantage to show clearly overpressure state and depression state. Indeed, $\rho_0$ being the volumetric mass of water, the overpressure state corresponds to  $\rho>\rho_0$ while $\rho<\rho_0$ represents a depression state.
The case $\rho=\rho_0$ is a critic one and a bifurcation point as we will see on \fig\ref{Depression}.

In the expression of the pressure (\ref{PressionLinearisee}), the sound speed is defined as
$\dsp c^2 = \frac{1}{\rho_0 \beta_0}$ where $\beta_0$ is the compressibility coefficient of water.
In practice, $\beta_0$ is $5.0\,10^{-10}\,m^2/N$ and thus $c\approx 1400\,m^2/s$.
 $p_a$ is some function  and without loss of generality, it may be set to zero.
Let us note that $p_a$ plays an important role in the construction of the mixed model \textbf{PFS} (as we will see in Section \ref{SectionCoupling}).

\noindent At the wet boundary, we assume a no-leak condition and we assume that the pipe is non-deformable. Thus we have the following crucial property:

\begin{equation}\label{MemeCourbePLane}
\left\{
\begin{array}{l}
\textrm{ If }(x,0,b_{\textrm{sl}}(x)) \textrm{ is the parametric representation of the plane curve }
\mathcal{C}_{\textrm{sl}} \textrm{ for }\\ \textrm{ free surface flows, }
\textrm{ then we define } \textrm{ continuously  the parametric} \\
\textrm{ representation} (x,0,b(x)) 
\textrm{ of the plane curve }\mathcal{C}_{\textrm{ch}} \textrm{ for pressurized flows.}
\end{array}
\right.
\end{equation}
As a consequence, the section $\Omega(x)$ (in pressurized state) orthogonal to the plane curve $\mathcal{C}_{\textrm{ch}}$ is a continuous extension of the free surface one. Henceforth, we note this curve $\mathcal{C}$.
At a given curvilinear abscissa, at the point $\omega\in\mathcal{C}$, we define pressurized section as follows:
$$\Omega(X) =
\left\{(Y,Z)\in\R^2; Z \in [-R(X), R(X)],\, Y \in \left[\alpha(X,Z),\beta(X,Z)\right]\right\}.$$
\begin{rque}
As the section is non-deformable,   $\Omega(x)$ depends only on the spatial variable $x$.
\end{rque}
Following the previous section, we proceed to the change of variables, namely we consider the application
$\mathcal{T} : (x,y,z)\rightarrow (X,Y,Z)$.

\subsection{Compressible Euler equations in curvilinear coordinates}\label{SubsectionEquationEulerCompressibleCurvilineaire}
\noindent Let  $(U,V,W)^t$ be the component of the fluid velocity  in variables $(X,Y,Z)$ given by $$(U,V,W)^t = \Uptheta (u,v,w)^t$$ where  $\Uptheta$ is the rotation matrix generated  around the axis $\ji$:
$$\Uptheta=\left(
\begin{array}{ccc}
 \cos\theta & 0 & \sin\theta \\
 0 & 1 & 0 \\
 -\sin\theta& 0 & \cos\theta \\
\end{array} \right)\,.$$

\subsubsection{Transformation of the equation of conservation of the mass}\label{SubsubsectionEquationEulerCompressibleMasseCurvilineaire}
Writing the  equation of conservation of the mass  (\ref{Euler3D_mass_conservation}) under a divergence form:
$$\dive_{t,x,y,z} \vecdeux{\rho}{\rho\UU}=0$$ and applying Lemma \ref{lemma}, we obviously get the equations in variables
$(X,Y,Z)$:
\begin{equation}\label{Mass}
\partial_t(J \rho ) +\partial_X (\rho U) + \partial_Y(J \rho V)+\partial_Z(J\rho W) = 0
\end{equation}
where $J$ is the determinant of the matrix $\mathcal{A}^{-1}$ (as already defined by \eqref{DefinitionA}).
\begin{rque}
Let us also remark that, from $(H)$, we   have $J(X,Z)>0$.
\end{rque}

\subsubsection{Transformation of the equation of conservation of the momentum}\label{SubsubsectionEquationEulerCompressibleMomentCurvilineaire}
Following Section \ref{SubsectionEulerCurvilinear}, namely:
\begin{itemize}
\item using Lemma \ref{lemma},
\item multiplying  the equation of conservation of the momentum (\ref{Euler3D_momentum_conservation}) on the left by the matrix $J \Uptheta$,
\end{itemize}
we get the equation for $U$ in the variables $(X,Y,Z)$:
\begin{equation}\label{Moment}
\partial_t(\rho J U) + \partial_X(\rho U^2) + \partial_Y(\rho JUV^2)+ \partial_Z(\rho JUW)+ \partial_X p  = -\rho
Jg\sin\theta(X)  +\rho U W  \dsp \frac{d\,}{dX}\cos\theta(X) .
\end{equation}
Other equations are unused since we want to derive a unidirectional model.
Let us also note, in the derivation of the \FS-model, all these equations were relevant to get the expression of the pressure.

\subsubsection{Formal asymptotic analysis}\label{SectionAnalysAsymptotiqueFormelleP}
As previously made in Section \ref{SectionAnalysAsymptotiqueFormelleSL}, we write the non-dimensioned version of equation equations  (\ref{Mass})-(\ref{Moment}) with respect to the small parameter $\varepsilon$ already introduced. In particular, we  seek for an approximation at main order with respect to the asymptotic expansion with respect to  $\varepsilon$. As pointed out before,  we will get $J\approx 1$ where $J$ is the determinant of the Jacobian matrix of the change of variables.

\noindent Let us recall that $$\dsp \frac{1}{\varepsilon} = \frac{L}{H}=
\frac{L}{l}= \frac{\overline{U}}{\overline{V}}= \frac{\overline{U}}{\overline{W}}$$ is assumed large enough where
$H$, $L$ and $l$ is a characteristic length of the height, the length,  and the width (where $l=H$ since we deal with only circular cross section pipe) and $\overline{U}$, $\overline{V}$ and
$\overline{W}$, is a characteristic speed following the main axis, the normal direction and the bi-normal one.
Then, let $T$ be a characteristic time such that $$\overline{U}=\frac{L}{T}.$$
We set the non-dimensioned  variables:
$$\widetilde{U} = \frac{U}{\overline U},\, \widetilde V = \varepsilon \frac{V}{ \overline U},\,
\widetilde W = \varepsilon \frac{W}{\overline U},\, $$ $$\widetilde X = \frac{X}{L},\, \widetilde Y =
\frac{Y}{H},\,
\widetilde Z = \frac{Z}{H},\,
\widetilde \rho = \frac{\rho}{\rho_0},\,\widetilde\theta = \theta.$$

\noindent With these notations, the non-dimensioned equations   (\ref{Mass})-(\ref{Moment}) are:
\begin{equation}\label{E3DC_X_Y_Z_rescaling}\left\{
\begin{array}{rcl}
\partial_{\widetilde t}(\widetilde\rho \widetilde J ) + \partial_{\widetilde X}({\widetilde\rho\widetilde U}) +
\partial_{\widetilde Y}(\widetilde\rho \widetilde J \widetilde V)+
\partial_{\widetilde Z}(\widetilde\rho \widetilde J \widetilde W)&=&0\\ \dsp
\partial_{\widetilde t}(\widetilde\rho \widetilde J \widetilde U) + \partial_{\widetilde
X}({\widetilde\rho\widetilde U}^2) +
\partial_{\widetilde Y}(\widetilde\rho \widetilde J \widetilde U \widetilde V)+
\partial_{\widetilde Z}(\widetilde\rho \widetilde J \widetilde U \widetilde W)+
\frac{1}{{M_a}^2}\partial_{\widetilde X} \widetilde\rho &=&\dsp
  \epsilon \widetilde\rho\widetilde U \widetilde W \widetilde \rho(\widetilde
  X)\\
& & \dsp -\widetilde\rho\displaystyle\frac{\sin\widetilde\theta(\widetilde X)}{{F_{r,L}}^2}\\
& & \dsp -  \frac{\widetilde Z}{{F_{r,H}}^2} \frac{d\,}{d\widetilde X}\cos\widetilde\theta(\widetilde X)
\end{array}\right.\end{equation}
where  $F_{r,M} = \dsp\frac{\bar U}{\sqrt{gM}}$ is the  Froude number and
$\dsp M_a = \frac{\overline{U}}{c}$ is the Mach number.

\noindent Formally, taking  $\varepsilon = 0 $ , the previous equations (\ref{E3DC_X_Y_Z_rescaling}) reads:
\begin{eqnarray}
\partial_{\widetilde t}(\widetilde \rho) + \partial_{\widetilde X}({\widetilde \rho \widetilde U}) +
\partial_{\widetilde Y}(\widetilde \rho \widetilde V)+
\partial_{\widetilde Z}(\widetilde \rho \widetilde W) &=& 0,  \label{E3DCCurvilinearRescaledEpsilonVanishes1}\\
\partial_{\widetilde  t}(\widetilde \rho \widetilde U  ) + \partial_{ \widetilde  X}(\widetilde \rho  {\widetilde
U}^2) + \partial_{ \widetilde Y}( \widetilde \rho  \widetilde U  \widetilde V)+
\partial_{\widetilde Z}( \widetilde \rho  \widetilde U   \widetilde W)+ \frac{1}{{M_a}^2} \partial_{  \widetilde X}
 \widetilde p
&= & -\widetilde\rho\displaystyle\frac{\sin\widetilde\theta(\widetilde X)}{{F_{r,L}}^2} \\
& & -   \frac{\widetilde Z}{{F_{r,H}}^2} \frac{d\,}{d\widetilde X}\cos\widetilde\theta(\widetilde X)
\label{E3DCCurvilinearRescaledEpsilonVanishes2}
\end{eqnarray}
\noindent\textcolor{red}{\danger} Henceforth, we note $(x,y,z)$ the dimensioned variables $(X,Y,Z)$, and $(u,v,w)$ dimensioned speed $(U,V,W)$.
Setting:
$$x=L \widetilde X,\,y=H \widetilde Y,\,z= H \widetilde Z,$$
$$u=\overline{U} \widetilde U,\,v=\epsilon \overline{U} \widetilde V ,\,v=\epsilon \overline{U} \widetilde W,$$
$\rho =  \, \widetilde \rho$, we finally get the following equations written in the curvilinear variables:
\begin{eqnarray}
\partial_t \rho + \partial_x (\rho u) + \partial_y (\rho v)+  \partial_z (\rho w)&=&0,
\label{E3DCCurvilinearRescaledDimensionnelles1} \\
\partial_t (\rho u) + \partial_x ( \rho u^2) + \partial_y (\rho u v)+ \partial_z ( \rho  u w)+ \partial_x p
&=&   -\dsp g\rho \sin\theta(x) -  g \rho z \frac{d\,}{dx}\cos\theta(x).
\label{E3DCCurvilinearRescaledDimensionnelles2}
\end{eqnarray}
where $p$ is the linearised pressure law given by (\ref{PressionLinearisee}).
\begin{rque}
The previous equations are the hydrostatic approximation of the Euler compressible equations where we have neglected the second and third momentum equation.
\end{rque}

\subsubsection{Vertical averaging of the hydrostatic approximation of Euler equations}
\label{SubsubsectionMoyennisationEulerCompressibleCurvilineaire}
The physical section of water, $S(x)$, and  the discharge, $Q(t,x)$, are defined by:
\begin{equation}\label{DefinitionS}
S(x) = \int_{\Omega(t,x)}\, dydz ,
\end{equation}
and
\begin{equation}\label{DefinitionQCh}
Q(t,x) = S(x) \overline{u}(t,x)
\end{equation}
where $\overline{u}$ is the mean speed over the section $\Omega(x)$:
\begin{equation}\label{Definitionu}
\overline{u}(t,x) = \frac{1}{S(t,x)}\int_{\Omega(t,x)} u(t,x,y,z)\, dydz .
\end{equation}
Let  $\mm \in \partial \Omega(x)$, we denote $\dsp \nn =
\frac{\mm}{|\mm|}$ the outward unit normal vector  to the boundary $\partial\Omega(x)$ at the point $\mm$ in the $\Omega$-plane and $\mm $ stands for the vector  \textbf{${\omega m}$} (c.f. \fig\ref{FigureOyOz}).

\noindent Following the section-averaging method performed to obtain the \FS-model, we integrate
System~(\ref{E3DCCurvilinearRescaledDimensionnelles1})-(\ref{E3DCCurvilinearRescaledDimensionnelles2})  over the cross-section $\Omega$.  Noting  the averaged values over $\Omega$ by the overlined letters  (except for $\overline{z}$), and using the approximations $\overline{\rho u}\approx \overline{\rho}\overline{u},\,\overline{\rho u^2}\approx
\overline{\rho}\overline{u}^2$, we get the following  shallow water like equations:
\begin{eqnarray}
\partial_{t}(\overline\rho S) + \partial_{x}({\overline\rho S \overline u}) & = & \dsp \int_{\partial \Omega(x)}
\rho
\left( u\partial_x \mm - \VV\right).\nn\, ds \label{STV_mass}\\
\partial_{t}(\overline\rho S \overline{u}) + \partial_x(\dsp\overline\rho S \overline{u}^2+c^2\overline\rho S) &=&
-\dsp{g\overline\rho S \sin \theta}+c^2\overline\rho \frac{d\, S}{dx} -   g\overline\rho S\overline z
\dsp \frac{d\,}{dx}\cos\theta\label{STV_momentum},\\
 &  & +   \dsp \int_{\partial \Omega(x)} \rho u \left(u\partial_x \mm - \VV\right).\nn\, ds\nonumber
\end{eqnarray}
where $\VV = (v,w)^t$ is the velocity field in the $(\NN, \BB)$-plane.
The integral terms appearing in (\ref{STV_mass}) and (\ref{STV_momentum}) vanish, as the pipe is infinitely rigid, i.e. $\Omega = \Omega(x)$ (see \cite{BG08} for more details about  deformable pipes). It follows the  non-penetration condition (see \fig\ref{Restriction}): $$\vectrois{u}{v}{w}.\nn_{\textbf{wb}} = 0 \,.$$

\noindent Omitting the overlined letters  (except for $\overline{z}$), setting  the conservative variables
\begin{eqnarray}
A = \dsp\frac{\rho }{\rho_0} S& \textrm{the FS \emph{equivalent wet area}} \label{EquivalentWetArea}\\
Q =\dsp A U & \textrm{ the FS \emph{equivalent discharge}} \label{EquivalentWetDischarge}\,.
\end{eqnarray}
and dividing Equations (\ref{STV_mass})-(\ref{STV_momentum}) by $\rho_0$, adding on each side of the equation (\ref{STV_momentum}) the quantity $-\dsp c^2\frac{d\, S}{dx}$, we get the pressurized model, called \Pres-model:
\begin{equation}\label{ModelP}
\left\{
\begin{array}{lll}
\dsp \partial_t A + \partial_{x} Q  &=&0, \\
\dsp \partial_t Q + \partial_x\left( \frac{Q^2}{A}+c^2 (A-S)\right) &=& -\dsp g A \sin\theta +c^2
\left(\frac{A}{S}-1\right) \frac{d\, S}{dx}  -   g A \overline z \dsp \frac{d\,}{dx}\cos\theta .
\end{array}
\right.
\end{equation}
\begin{rque}
In term of area,  a depression occurs when $A<S$ (i.e $\rho<\rho_0$) and an overpressure if $A>S$ (i.e. $\rho>\rho_0$).
\end{rque}
As introduced previously for the \textbf{FS}-model, we may introduce the friction term  $-\rho g S_f\, \TT$ given by the Manning-Strickler law (see e.g. \cite{SWB98}): $$S_f(S,U)=K(S)U|U|$$ where $K(S)$ is defined by: $\dsp K(S) = \frac{1}{K_s^2 R_h(S)^{4/3}}$, $K_s>0$ is the Strickler coefficient
of roughness depending on the material and $R_h(S)= S/P_m$ is the hydraulic radius where $P_m$ is the perimeter of the wet surface area (length of the part of the channel's section in contact with the water, equal to $2\,\pi\,R$ in the case of circular pipe).

\section{The \PFS-model}\label{SectionCoupling}
In the previous sections, we have proposed two models, one for the free surface flows and the other for pressurized (compressible) flows which are very close to each other. In this section,  we are motivated to connect ``continuously" these two models through the transition points.
 Let us recall these two models.
\paragraph*{The \FS-model.\newline}
\begin{equation}\label{freesurfacesl}
\left\{
\begin{array}{rcl}
\dsp \partial_t A_{sl} + \partial_x Q_{sl} &=&0,\\
\dsp \partial_t Q_{sl} + \partial_x \left(\dsp\frac{Q_{sl}^2}{A_{sl}}+p_{sl}(x,A_{sl})\right) &=& -\dsp g A_{sl}\dsp
\frac{d\, Z}{dx}  + Pr_{sl}(x,A_{sl}) -g A_{sl} \overline z \dsp \frac{d\,}{dx}\cos\theta\\
& & \dsp  - g K(x,A_{sl})\dsp\frac{Q_{sl}|Q_{sl}|}{A_{sl}}
\end{array}\right.
\end{equation}
where
$$
\begin{array}{lll}
p_{sl}(x,A_{sl}) &=& \dsp g I_1(x,A_{sl}) \cos\theta, \\
Pr_{sl}(x,A_{sl}) &=& \dsp g I_2(x,A_{sl}) \cos\theta, \\
K(x,A_{sl}) &=& \dsp \frac{1}{K_s^{2} R_h(A_{sl})^{4/3}} \ ,\\
\end{array}
$$
with $I_1$ and $I_2$ are defined by (\ref{DefinitionI1}) and (\ref{DefinitionI2}).
\paragraph*{The \Pres-model.\newline}
\begin{equation}\label{pressurisedch}
\left\{
\begin{array}{rcl}
\dsp \partial_t A_{ch} + \partial_x Q_{ch} &=&0,\\
\dsp \partial_t Q_{ch} + \partial_x \left(\dsp\frac{Q_{ch}^2}{A_{ch}}+p_{ch}(x,A_{ch})\right) &=& -\dsp g A_{ch}\dsp
\frac{d\, Z}{dx}  + Pr_{ch}(x,A_{ch}) -g A_{ch} \overline z \dsp \frac{d\,}{dx}\cos\theta\\
& & \dsp  - g K(x,S)\dsp\frac{Q_{ch}|Q_{ch}|}{A_{ch}}
\end{array}\right.
\end{equation}
where
$$
\begin{array}{lll}
p_{ch}(x,A_{ch})  &=& \dsp c^2(A_{ch}-S), \\
Pr_{ch}(x,A_{ch}) &=& \dsp c^2 \left(\frac{A_{ch}}{S}-1\right)\frac{d\, S}{dx} , \\
K(x,S) &=& \dsp \frac{1}{K_s^{2} R_h(S)^{4/3}} \ .\\
\end{array}
$$

We remark that the term $\dsp\frac{d\, Z}{dx}$, $\overline{z}\frac{d}{dx}\cos\theta$ and the friction are similar in both models (where we set $Z(x)=b(x)$).
\begin{rque}\label{RemarkRestrictionPlaneCurve}

\noindent The plane curve with parametrization $(x,0,b(x))$ is chosen as the main pipe axis in the axis-symmetric case. Actually this choice is the more convenient for pressurised flows while the
bottom line is adapted to free surface flows. Thus, we must assume small variations of the section
($\dsp S'$ small) or equivalently small angle $\varphi$ as displayed on \fig\ref{Restriction} in order to keep a continuous connection of the term $Z(x)$ from one to other type of flows.
\begin{figure}[H]
 \begin{center}
 \includegraphics[scale = 0.6]{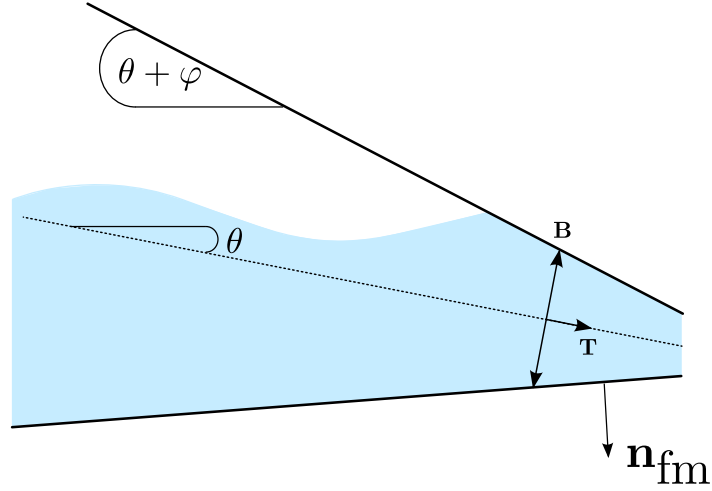}
  \caption{Some restriction concerning the geometric domain: $\varphi < \theta$.}
  \label{Restriction}
 \end{center}
\end{figure}
\end{rque}
The real difference between these two models is mainly due to the pressure law: one is  of ``acoustic" type while the other is hydrostatic. We are then motivated by defining a suitable couple of ``mixed" variables in order to connect ``continuously" these two models through the transition points. But, necessarily, the gradient of flux of the new system will be discontinuous at transition points, due to the difference of sound speed (as we will see below).

\begin{figure}[H]
\centering
\includegraphics[scale=0.55]{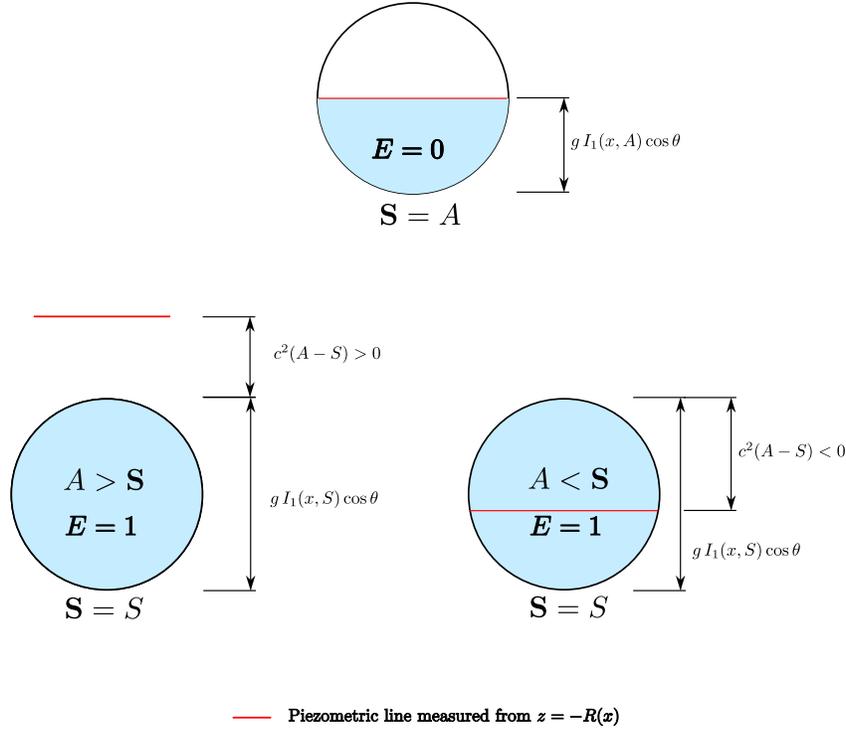}
\caption{Free surface state $p(X,A,0) = g\,I_1(X,A)\cos\theta$ (top), pressurized state with overpressure
$p(x,A,1)>0$ (bottom left), pressurized state with depression $p(x,A,1)<0$ (top right).}\label{Depression}
\end{figure}

\noindent To this end, we introduce a state indicator $E$ (see \fig\ref{Depression}) such that:
\begin{equation}\label{DefinitionE}
E =\left\{
\begin{array}{ll}
1 & \textrm{ if the state is pressurised: } (\rho\neq\rho_0), \\
0 & \textrm{ if the state is free surface: } (\rho=\rho_0).\\
\end{array}
\right.
\end{equation}

\begin{figure}[H]
\centering
\includegraphics[scale=0.55]{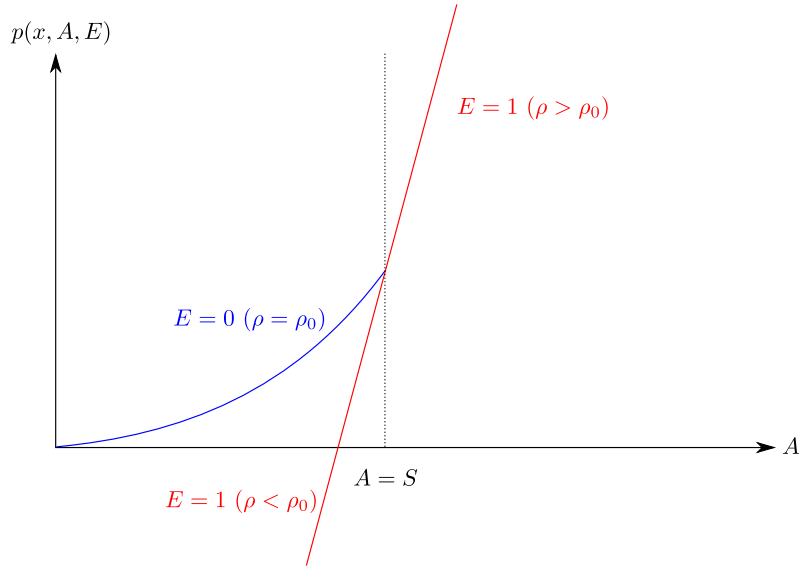}
\caption{Pressure law in the case of pipe with circular section.}\label{MultivaluedPressure}
\end{figure}

\noindent Next, we define the \emph{physical wet area} $\SE$ by:
\begin{equation}\label{DefinitionPhysicalWetArea}
\SE = \SE(A,E) = \left\{
\begin{array}{lll}
S       & \textrm{ if } & E = 1, \\
A       & \textrm{ if } & E = 0.
\end{array}
\right.
\end{equation}
We introduce then a couple of variables, called ``mixed variables":
\begin{eqnarray}
A = \dsp\frac{\rho}{\rho_0} \SE, \label{EquivalentWetAreaPFS}\\
Q =\dsp A u\label{EquivalentWetDischargePFS}
\end{eqnarray}
which satisfies:
\begin{itemize}
 \item if the flow is free surface, $\rho=\rho_0$, $E=0$ and consequently $\SE=A$, and
\item  if the flow is pressurized, $\rho\neq\rho_0$, $E=1$ and consequently  $\SE=S$.
\end{itemize}
\noindent To construct a ``mixed" pressure law (c.f. \fig\ref{Depression}), we set:
\begin{equation}\label{PFSModelPressureInFlux}
 p(x,A,E) =  c^2(A-\SE) + gI_1(x,\SE)\cos\theta
\end{equation}
where the term $I_1$ is given by:
\begin{equation}\label{DefinitionI1PFS}
I_1(x,\SE) = \dsp\int_{-R(x)}^{\mathcal{H}(\SE)}
(\mathcal{H}(\SE)-z)\, \sigma(x,z)\,dz
\end{equation}
with $\mathcal{H}$ representing   the $z$-coordinate of the water level:
\begin{equation}\label{DefinitionPhysicalWaterHeight}
\mathcal{H} = \mathcal{H}(\SE)=\left\{
\begin{array}{lll}
h(A)       & \textrm{ si } & E = 0, \\
R(x)       & \textrm{ si } & E = 1.
\end{array}
\right.
\end{equation}
Thus, the constructed pressure is continuous throughout the transition points:
$$\lim_{\stackrel{A\to S}{A<S}} p(x,A,E) = \lim_{\stackrel{A\to S}{A>S}} p(x,A,E)$$
but its gradient is discontinuous (c.f. \fig\ref{MultivaluedPressure}):
$$\frac{\partial p}{\partial A}(x,A,0) = \sqrt{\frac{g A}{T}} \neq c^2=\frac{\partial p}{\partial A}(x,A,1).$$
\begin{rque}
The transition point ($\rho=\rho_0$) is then a bifurcation point.
\end{rque}

\noindent From the \FS-model \eqref{freesurfacesl}, the \Pres-model \eqref{pressurisedch}, the ``mixed"  variables (\ref{EquivalentWetAreaPFS})-(\ref{EquivalentWetDischargePFS}), the state indicator $E$ (\ref{DefinitionE}), the physical height of water $\SE$ \eqref{DefinitionPhysicalWetArea} and  the pressure law (\ref{PFSModelPressureInFlux}), we can define the \PFS-model (\textbf{P}ressurised and \textbf{F}ree \textbf{S}urface) for unsteady mixed flows in closed water pipes with variable section and slope, as follows:
\begin{equation}\label{PFS}
\left\{
\begin{array}{rcl}
\partial_t A + \partial_x Q &=&0\\
\partial_t Q + \partial_x\left(\dsp\frac{Q^2}{A}+p(x,A,E)\right) &=&
-\dsp g A \frac{d\, Z}{dx}  + Pr(x,A,E)\\ & & -G(x,A,E)\\ & &- g K(x,A,E)\dsp\frac{Q|Q|}{A}
\end{array}\right.
\end{equation}
\noindent where  $Pr$, $K$, and $G$ represent respectively the pressure source term, the curvature term and the friction:
$$\begin{array}{lll}
Pr(x,A,E) &=& \dsp c^2\left(\frac{A}{\SE}-1\right)\, \frac{d\, S}{dx}   + g I_2(x,\SE) \cos\theta \\
          &\textrm{ with } &  I_2(x,\SE) = \dsp\int_{-R(x)}^{\mathcal{H}(\SE)}
(\mathcal{H}(\SE)-z)\,\partial_x \sigma(x,z)\,dz,\\
 G(x,A,E) &=& \dsp g A\,  \overline{z}(x,\SE) \dsp \frac{d\,}{dx}\cos\theta, \\
K(x,A,E) &=& \dsp \frac{1}{K_s^{2} R_h(\SE)^{4/3}} .
  \end{array}
$$
\begin{rque}\label{RqueRecoveringBothModels}
From the \textbf{PFS} equations, it is easy to recover the \FS-model. Indeed, if $E=0$ then
$\SE(A,E) = A$ and the pressure law (\ref{PFSModelPressureInFlux}) is the hydrostatic pressure, the term $Pr$ is also the
hydrostatic pressure source term. It follows that the \PFS-system (\ref{PFS}) coincide exactly with the \FS-model (\ref{freesurfacesl}).

\noindent If $E=1$, then $\SE(A,E) = S$ and the pressure law gives $\dsp c^2(A-S) + gI_1(x,S)\cos\theta$ which is exactly the linearised pressure law \emph{if we consider $p_a(x) = gI_1(x,S)\cos\theta$ instead of $0$}.
We can then see the term $p_a$ as a limit state between an over-pressurised zone  and a de-pressurized one.
We show different situations on \fig\ref{Depression}.
\end{rque}
\begin{rque}
We have seen that when the flow is fully pressurized, the overpressure states are reported when $A>S$ and depression states when $ A <S $. However, when the flow is mixed and $A<S $, we do not know \emph{a priori} if the state is free surface or pressurized. Therefore,  the indicator state $E$ is there to overcome this difficulty.
Thus, combining this with a discrete algorithm on $E$ is useful to describe both depression areas and free surface ones.
When $A>S$, without any ambiguity, the pressurised state is proclaimed.
From a numerical point of view, the transition points between two types of flows are treated as a free boundary, corresponding to the discontinuity of the gradient  of the pressure (for more details, see \cite{BEG09_2,TheseErsoy}).
\end{rque}

\noindent The \textbf{PFS}-model (\ref{PFS}) satisfies the following properties:
\begin{thm}\label{ThmPFSModel}
\begin{enumerate}
\item[]
\item The right eigenvalues of System (\ref{PFS}) are given by:
$$\lambda^- =u-c(A,E),\,\lambda^+ = u + c(A,E)
$$
with $c(A,E)= \left\{
\begin{array}{lll}
\dsp\sqrt{g\,\frac{A}{T(A)}\,\cos\theta} &\textrm{ if } & E=0,\\
\dsp c  &\textrm{ if } & E=1 .
\end{array}
\right.
$

\noindent Then, System (\ref{PFS}) is strictly hyperbolic on the set:
$$\left\{A(t,x)>0\right\}\,.$$
\item For smooth solutions, the mean velocity $u = Q/A$ satisfies
\begin{equation}\label{ThmPFSEquationForU}
\begin{array}{c}
\partial_t u + \partial_x \left(\dsp\frac{ u^2}{2} + c^2 \ln(A/\SE) + g\mathcal{H}(\SE)\cos\theta + gZ \right) \\ =
-g K(x,A,E) u|u| \leqslant 0.
\end{array}
\end{equation}
\noindent The quantity $\dsp \frac{ u^2}{2} + c^2 \ln(A/\SE) + g\mathcal{H}(\SE)\cos\theta + gZ$ is called the
total head.
\item The still water steady state reads:
\begin{equation}\label{ThmPFSSteadyState}
u = 0 \;\mbox{ and } \; c^2 \ln(A/\SE) + g\mathcal{H}(\SE)\cos\theta + g Z = cte.
\end{equation}
\item It admits a mathematical entropy
\begin{equation}\label{ThmPFSMathematicalEntropy}
\mathcal{E}(A,Q,E) =\dsp \frac{Q^2}{2A} + c^2 A \ln(A/\SE)+ c^2 S + g A
\overline{z}(x,\SE)\cos\theta + gA Z
\end{equation}
\noindent which satisfies the entropy relation for smooth solutions
\begin{equation}\label{ThmPFSEntropy}
\partial_t \mathcal{E} +\partial_x \Big(\left(\mathcal{E}+p(x,A,E)\right)u\Big) = -g A K(X,A,E) u^2 |u| \leqslant 0\,.
\end{equation}
\end{enumerate}
\end{thm}
Notice that the total head and $\mathcal{E}$ are defined continuously through the transition points.

\begin{rque}
The  term $A\overline{z}(x,A) (\cos\theta)'$ is also called ``corrective term'' since it allows to write the  equation
(\ref{ThmPFSEquationForU}) and (\ref{ThmPFSEntropy}) with (\ref{ThmPFSMathematicalEntropy}).
\end{rque}

\noindent \textbf{Proof of Theorem~\ref{ThmPFSModel}:} The results~(\ref{ThmPFSEquationForU}) and (\ref{ThmPFSEntropy})  are obtained in a
classical way. Indeed, Equation (\ref{ThmPFSEquationForU}) is obtained by subtracting the result of the multiplication of the mass  equation
by $u$ to the momentum equation. Then multiplying the mass  equation by $\left(\dsp\frac{ u^2}{2} + c^2 \ln(A/\SE) +
g\mathcal{H}(\SE)\cos\theta + gZ \right)$ and adding the result of the multiplication of Equation (\ref{ThmPFSEquationForU}) by $Q$,
we get:
$$\begin{array}{l}
\partial_t \left(\dsp \frac{Q^2}{2A} +c^2 A \ln(A/\SE)+ c^2 S +g A \overline{z}(x,\SE)\cos\theta + gAZ\right) \\
+\partial_x
\left(\left(\dsp \frac{Q^2}{2A} + c^2 A \ln(A/\SE)+ c^2 S +gA \overline{z}(x,\SE)\cos\theta +
gAZ+p(x,A,E)\right)u\right) \\
+c^2\left(\dsp\frac{A}{\SE}-1\right)\partial_t \SE = -gAK(x,A,E) u^2 |u| \leqslant 0 \,.\end{array}$$ We see that
the
term  $c^2\left(\dsp\frac{A}{\SE}-1\right)\partial_t \SE$ is identically $0$ since we have  $\SE=A$   when
the flow is free
surface whereas $\SE=S(x)$ when the flow is pressurised. Moreover, from the last inequality, when $\SE=A$, we have the
classical entropy
inequality (see \cite{BG07,BG08}) with $\mathcal{E}$: $$\mathcal{E}(A,Q,E) =\dsp \frac{Q^2}{2A} +gA \overline{z}(x,A)\cos\theta +
gAZ $$
while in the pressurised case, it is: $$\mathcal{E}(A,Q,E) =\dsp \frac{Q^2}{2A} + c^2 A \ln(A/S)+ c^2 S  + gAZ.$$
Finally, the entropy for the \textbf{PFS}-model reads: $$\mathcal{E}(A,Q,E) =\dsp \frac{Q^2}{2A} + c^2 A \ln(A/\SE)+ c^2 S  +gA
\overline{z}(X,\SE)\cos\theta + g A  Z.$$
Let us remark that the term $c^2 S$ makes $\mathcal{E}$ continuous through transition points and it permits also to write   the
entropy flux under the classical form $(\mathcal{E}+p)u$.
 \begin{flushright}
  $\blacksquare$
 \end{flushright}

\section{Perspectives}
In view of the difference of sound speed ($c\approx 1400\:m/s$ for a pressurised state and $c\approx 1\:m/s$ for a free surface state), the gradient of the pressure, thus the flux of the \PFS~equations, is discontinuous throughout the transition points. More generally, these equations belong to a class of hyperbolic systems of  conservation laws with discontinuous gradient, especially a generalization of equations coupled through a fixed discontinuity (see \cite{LWR55, LWR56, Mochon87} with the classical example of Lighthill-Whitham-Richards model for road traffic) since, in the present case the discontinuity is mobile. It is then an interesting and a difficult problem because of the definition of the solution
associated to the Riemann problem. In general, given two initial states which are not connected by a shock wave, there exits an infinite number of paths to connect them through the interface.
For instance, in Boutin's thesis \cite{Boutin09}, he defined paths using physical criteria that enable to extract the solution. To our knowledge, and up to  date, there are no results for the mobile discontinuities and the \PFS-model is a nice example of such open problem.
However, in each region where the gradient of the flux is continuous,  the solution is constructed in a classical
way (see, for example, \cite{Toro92}).

\section*{Acknowledgements}
This work is supported by the ``Agence Nationale de la Recherche'' referenced by  ANR-08-BLAN-0301-01 and 
the  second author was  supported by  the ERC Advanced Grant FP7-246775 NUMERIWAVES.
\bibliographystyle{plain}

\end{document}